\documentclass[oneside]{amsart}
\newcommand{\formatswitch}{preprint}

\usepackage{amssymb}
\usepackage{ifthen}
\usepackage{verbatim}
\usepackage{psfig}
\usepackage[all]{xypic}

\newcommand{\tref}[1]{(\ref{#1})}

\ifthenelse{\equal{\formatswitch}{big}}{
\setlength{\textwidth}{432pt}
\setlength{\oddsidemargin}{18.8775pt}

\setcounter{tocdepth}{1}
}{
\ifthenelse{\equal{\formatswitch}{paper}}{

\setcounter{tocdepth}{1}
}{
\setcounter{tocdepth}{2}
}
}
\DeclareMathAlphabet\EuScript{U}{eus}{m}{n}
\DeclareMathAlphabet\EuScriptb{U}{eus}{b}{n}

\newcommand{\claimenum}{\renewcommand{\theenumi}{\alph{enumi}}
 \renewcommand{\labelenumi}{\textit{(\theenumi)}}
 \renewcommand{\theenumii}{\roman{enumii}}
 \renewcommand{\labelenumii}{\textit{(\theenumii)}}
 \begin{enumerate}}
\newcommand{\claimenumend}{\end{enumerate}}

\newcommand{\romanenum}{\renewcommand{\theenumi}{\roman{enumi}}
 \renewcommand{\labelenumi}{\textit{(\theenumi)}}
 \renewcommand{\theenumii}{\alph{enumii}}
 \renewcommand{\labelenumii}{\textit{(\theenumii)}}
 \begin{enumerate}}
\newcommand{\romanenumend}{\end{enumerate}}

\newtheorem{dummy}{realdumb}[section]
\newtheorem{thm}{Theorem}

\newtheorem{lemma}[dummy]{Lemma}
\newtheorem{prop}[dummy]{Proposition}
{\theoremstyle{definition} }

{\theoremstyle{definition} }
\newtheorem{cor}{Corollary}[dummy]

\renewcommand{\text}{\mathrm}

\newcommand{\strutdepth}{\dp\strutbox}
\newcommand{\marginalnote}[1]
   {\strut\vadjust{\kern-\strutdepth\domarginalnote{#1}}}
\newcommand{\domarginalnote}[1]{\vtop to \strutdepth{
  \baselineskip\strutdepth
   \vss\llap{ #1\ \ }\null}}  %get sevenpoint font here
\newcounter{showlabelflag}
\newcounter{makelabelflag}
\newcommand{\showlabels}{\setcounter{showlabelflag}{1}}

\newcommand{\makelabels}{\setcounter{makelabelflag}{1}}
\newcommand{\hidelabels}{\setcounter{showlabelflag}{2}}
\newcommand{\mylabel}[1]{
  \ifthenelse{\value{makelabelflag}=1}
    {\label{#1}}{}
  \ifthenelse{\value{showlabelflag}=1}
    {\marginpar{#1}}{}\relax}

\newcommand{\R}{{\mathbf R}}

\newcommand{\Z}{{\mathbf Z}}
\newcommand{\N}{{\mathbf N}}
\newcommand{\sub}{\subseteq}

\ifthenelse{\equal{\formatswitch}{draft}}{
\showlabels
}{\relax}

\newcommand{\opi}{\overline{\pi}}

\newcommand{\supp}{\mathrm{Supp}}

\newcommand{\mymargin}[1]{
  \ifthenelse{\value{showlabelflag}=1}
    {\marginpar{#1}}{}\relax}

\newcounter{enumo}

\newcommand{\RRsh}{\kern -1 pt \Rsh}

\newcounter{keepitemnum}

\newcounter{keepitemnumm}

\begin{document}

\bibliographystyle{amsplain}
\begin{center}{\bfseries Higher Dimensional Thompson Groups\footnote{AMS
Classification (2000): primary 20B27, secondary 20E32, 20F55, 
57S25, 20B22, 37E05, 37E15.}}\end{center}
\vspace{3pt}
\begin{center}{MATTHEW G. BRIN}\end{center}
\vspace{4pt}
\vspace{3pt}
%\begin{center}\today\end{center}
\begin{center}January 2, 2004\end{center}

\tableofcontents
\CompileMatrices

%\makelabels has to come AFTER \tableofcontents 
%to prevent the labels in the section headers from 
%being processed twice

\makelabels
%\showlabels
\hidelabels

%**end of header

\section{Introduction}

Three groups \(F\subseteq T\subseteq V\), known as Thompson groups,
have generated interest since R.~J.~Thompson introduced them in the
late 1960s.  Part of their initial interest was the fact that \(T\)
and \(V\) supplied the first known examples of infinite, simple,
finitely presented groups.  Since then, other properties of the
groups have been studied as well as their interaction with other
areas of mathematics.  The standard reference for the groups
\(F\subseteq T\subseteq V\) is \cite{CFP}.  It is not necessary to
have \cite{CFP} in hand while reading this paper, but it would not
hurt.

The largest group \(V\) can be described as a subgroup of the
homeomorphism group of the Cantor set \(C\).  Intrinsic to this
description is the standard ``deleted middle thirds'' construction
of the Cantor set as a subset of the unit interval.  Since the unit
interval is a 1-dimensional object, we will refer to \(V\) as a
1-dimensional Thompson group.

In this paper, we describe an intrinsically 2-dimensional group
\(2V\) that is more naturally described as a subgroup of the
homeomorphism group of \(C \times C\).  We will show in this paper
that \(2V\) is infinite, simple and finitely generated, and in
another paper \cite{brin:hd4}, that \(2V\) is finitely presented.

In spite of the fact that \(C\times C\) is homeomorphic to \(C\), we
will also show in this paper that \(2V\) is not isomorphic to \(V\).
In fact, we will show that \(2V\) is not isomorphic to any group in
a list of other infinite, simple, finitely presented groups that are
closely related to \(T\) and \(V\).  However, we do not show that
\(2V\) is not isomorphic to all known infinite, simple, finitely
presented groups.

The group \(2V\) is a member of a family of groups \(nV\) of which
\(1V=V\).  We describe these groups for \(n>2\), but prove little
about them.  They are all infinite and it is a reasonable guess that
they are also simple and finitely presented.  Part of the charm of
the family of groups \(nV\) is the uncanny resemblance they have to
the little cubes operads of \cite{may:geom+loop} and
\cite{boardman+vogt}.  The importance of this resemblance is
speculative.

It would be nice to claim that \(2V\) is a \emph{new} infinite,
simple, finitely presented group.  As the list of such groups grows
longer, such claims become less interesting and harder to support.
The difficulty of establishing such a claim is increased by the lack
of an official list.  So we make no such claim and introduce the
groups here for their pleasant properties and their potential
applicability to other areas of mathematics.  See the discussion in
Sections \ref{RubinScopeSec} and \ref{NoIsoScopeSec} for references
to all infinite, simple, finitely presented groups that the author
is aware of.

The proof that \(2V\) is not isomorphic to \(V\) and other simple
groups is not algebraic.  We invoke a difficult theorem of M. Rubin
that gives conditions under which isomorphism implies topological
conjugacy.  Then we demonstrate irreconcilable differences between
the dynamics of \(V\) and \(2V\): there is an element of \(2V\)
exhibiting ``chaotic'' behavior and there is no such element of
\(V\).  Thus part of this paper studies properties of the older
group \(V\).  Dynamical properties of Thompson groups have been
studied before in \cite{ghys+serg} and \cite{brin:cham}.

The chaotic element in \(2V\) is a well known map often referred to
as the baker's map.  Some aspects of the baker's map are discussed
briefly in the last two sections of the paper.  In particular we
show that the baker's map is a product of commutators of some rather
easily described homeomorphisms of the Cantor set.

The paper \cite{CFP} (in Sections 1 and 2 and the first few pages of
5 and 6) covers the representation of the groups \(F\), \(T\) and
\(V\) as homeomorphism groups and also the encoding of the elements
of the groups by pairs of trees.  We will need the pairs of trees
point of view in analyzing the dynamics of elements of \(V\).

\section{Defining \protect\(2V\protect\)}\mylabel{DefTwoVeeSec}

We will define the group \(2V\) as a subgroup of the group
\(H(C^2)\) of self homeomorphisms of \(C^2=C\times C\).  We
introduce some conventions that we will use.

In this paper, homeomorphisms will act on the left and will be
composed from right to left.  The \emph{support} of an element \(h\)
of \(2V\), denoted \(\supp(h)\), will be the closure of \(\{x\in
C\times C\mid x\ne h(x)\}\).  We will use \(h^g\) to denote
\(ghg^{-1}\) and from this, we will have \(\supp(h^g) =
g\left(\supp(h)\right)\).  We will use \([h,g]\) to denote the
commutator \(\displaystyle{hgh^{-1}g^{-1} =
h\left(h^{-1}\right)^g}\).

We will think of the Cantor set \(C\) as the usual deleted middle
thirds subset of the unit interval \(I=[0,1]\) and we will identify
elements of \(C\) with infinite words in the alphabet \(\{0,1\}\).
These words are just functions from the natural numbers
\(\N=\{0,1,2,\cdots\}\) to \(\{0,1\}\).  The identification of these
words with elements of \(C\) takes words beginning with 0 to points
in \([0,\frac13]\), words beginning with 1 to points in
\([\frac23,1]\), words beginning with 00 to points in
\([0,\frac19]\) and so forth.  If \(w\) is an infinite word in
\(\{0,1\}\), then doubling all the ``digits'' in \(w\) gives the
ternary expansion as a real number of the image of \(w\) in \(C\).

We view \(C\times C\) as a subset of the unit square
\(I^2=[0,1]\times [0,1]\).  We next develop methods of describing
individual elements of \(2V\).

\subsection{Patterns of rectangles}

We will describe certain partitions of \(C\times C\) by patterns of
rectangles in the unit square.  First we describe the patterns, and
then we say what partitions of \(C\times C\) the patterns are
describing.

We inductively define what we call \emph{patterns} in \(I^2\).  Each
pattern is a finite set of rectangles in \(I^2\) with pairwise
disjoint, non-empty interiors, with sides parallel to the sides of
\(I^2\) and whose union is all of \(I^2\).  The first pattern is
\(I^2\) itself and will be called the \emph{trivial pattern}.  If
\(P\) is a pattern and \(R\) is one of the rectangles in \(P\), then
we can create a new pattern by removing \(R\) from \(P\) and
replacing it with two congruent rectangles \(R'\) and \(R''\)
obtained from \(R\) by dividing \(R\) exactly in half by either a
horizontal line, or a vertical line.  A pattern is one that can be
obtained from the trivial pattern in this way in a finite number of
steps.

Below we show four patterns.  Pattern \((a)\) is the trivial pattern
and patterns \((c)\) and \((d)\) are obtained from pattern \((b)\)
by the two possible ways of dividing the upper left rectangle in
\((b)\).

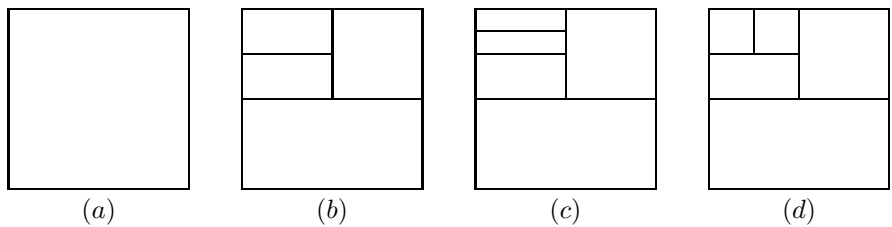
\[
\xy (0,0); (0,24)**@{-}; (24,24)**@{-}; (24,0)**@{-};
(0,0)**@{-};
(12,-3)*{(a)};
\endxy
\qquad
\xy (0,0); (0,24)**@{-}; (24,24)**@{-}; (24,0)**@{-};
(0,0)**@{-};
(0,12); (24,12)**@{-}; (12,12); (12,24)**@{-};
(0,18); (12,18)**@{-};
(12,-3)*{(b)};
\endxy
\qquad
\xy (0,0); (0,24)**@{-}; (24,24)**@{-}; (24,0)**@{-};
(0,0)**@{-};
(0,12); (24,12)**@{-}; (12,12); (12,24)**@{-};
(0,18); (12,18)**@{-};
(0,21); (12,21)**@{-};
(12,-3)*{(c)};
\endxy
\qquad
\xy (0,0); (0,24)**@{-}; (24,24)**@{-}; (24,0)**@{-};
(0,0)**@{-};
(0,12); (24,12)**@{-}; (12,12); (12,24)**@{-};
(0,18); (12,18)**@{-};
(6,18); (6,24)**@{-};
(12,-3)*{(d)};
\endxy
\]

We note that different sequences of modifications to the trivial
pattern can result in the same pattern.  Several of the fundamental
relations that we will work with come from the several sequences of
operations that yield the pattern below when starting from the
trivial pattern.

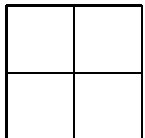
\[
\xy
(-9,-9); (-9,9)**@{-}; (9,9)**@{-}; (9,-9)**@{-}; (-9,-9)**@{-};
(-9,0); (9,0)**@{-}; (0,-9); (0,9)**@{-};
\endxy
\]

We next number our patterns.  A \emph{numbered pattern} is a
pattern with a one-to-one correspondence between \(\{0,1,\cdots,
n-1\}\) and the rectangles in the pattern where \(n\) is the number
of rectangles in the pattern.  Below we show two different numbered
patterns based on the same pattern.

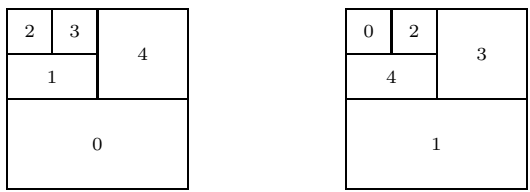
\[
\xy (0,0); (0,24)**@{-}; (24,24)**@{-}; (24,0)**@{-};
(0,0)**@{-};
(0,12); (24,12)**@{-}; (12,12); (12,24)**@{-};
(0,18); (12,18)**@{-};
(6,18); (6,24)**@{-};
(12,6)*{\scriptstyle0};
(6,15)*{\scriptstyle1};
(3,21)*{\scriptstyle2};
(9,21)*{\scriptstyle3};
(18,18)*{\scriptstyle4};
\endxy
\qquad
\qquad
\qquad
\xy (0,0); (0,24)**@{-}; (24,24)**@{-}; (24,0)**@{-};
(0,0)**@{-};
(0,12); (24,12)**@{-}; (12,12); (12,24)**@{-};
(0,18); (12,18)**@{-};
(6,18); (6,24)**@{-};
(12,6)*{\scriptstyle1};
(6,15)*{\scriptstyle4};
(3,21)*{\scriptstyle0};
(9,21)*{\scriptstyle2};
(18,18)*{\scriptstyle3};
\endxy
\]

\subsection{Partitions from patterns}

Each rectangle in a pattern corresponds to a closed and open subset
of \(C\times C\).  We describe this correspondence inductively.

The rectangle \(I^2\) corresponds to all of \(C\times C\).

Let a rectangle \(R\) in a pattern correspond to the subset \(A(R)\)
in \(C\times C\).  Let \(R'\) and \(R''\) be obtained from \(R\) by
dividing \(R\) equally by a vertical line, with \(R'\) the left
rectangle and \(R''\) the right.  Then \(R'\) corresponds to the
left third of \(A(R)\) and \(R''\) corresponds to the right third of
\(A(R)\).

Similarly, if \(R'\) and \(R''\) be obtained from \(R\) by dividing
\(R\) equally by a horizontal line, with \(R'\) the bottom rectangle
and \(R''\) the top.  Then \(R'\) corresponds to the bottom third of
\(A(R)\) and \(R''\) corresponds to the top third of \(A(R)\).

The next figures give the basic correspondences for the patterns
obtainable from the trivial pattern \(I^2\) using one division.

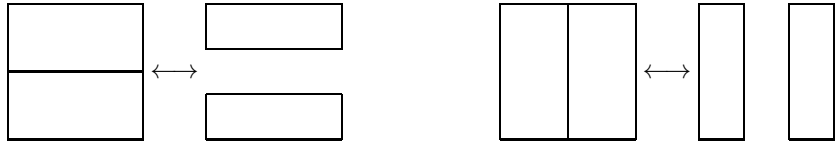
\[
\xy
(-9,-9); (-9,9)**@{-}; (9,9)**@{-}; (9,-9)**@{-}; (-9,-9)**@{-};
(-9,0); (9,0)**@{-};
\endxy
\longleftrightarrow
\xy
(-9,-9); (-9,-3)**@{-}; (9,-3)**@{-}; (9,-9)**@{-}; (-9,-9)**@{-};
(-9,9); (-9,3)**@{-}; (9,3)**@{-}; (9,9)**@{-}; (-9,9)**@{-};
\endxy
\qquad
\qquad
\qquad
\xy
(-9,-9); (-9,9)**@{-}; (9,9)**@{-}; (9,-9)**@{-}; (-9,-9)**@{-};
(0,-9); (0,9)**@{-};
\endxy
\longleftrightarrow
\xy
(-9,-9); (-9,9)**@{-}; (-3,9)**@{-}; (-3,-9)**@{-}; (-9,-9)**@{-};
(3,-9); (3,9)**@{-}; (9,9)**@{-}; (9,-9)**@{-}; (3,-9)**@{-};
\endxy
\]

Thus it is seen that a pattern is just a lazy way of drawing a
particular partition of \(C\times C\) into closed and open sets in
\(C\times C\).  Each rectangle in the pattern gives one set in the
partition.  An example is pictured below.

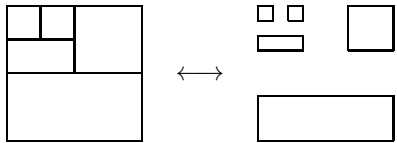
\[
\xy (-9,-9); (-9,9)**@{-}; (9,9)**@{-}; (9,-9)**@{-};
(-9,-9)**@{-};
(-9,0); (9,0)**@{-}; (0,0); (0,9)**@{-};
(-9,4.5); (0,4.5)**@{-};
(-4.5,4.5); (-4.5,9)**@{-};
\endxy
\quad
\longleftrightarrow
\quad
\xy
(-9,-9); (-9,-3)**@{-}; (9,-3)**@{-}; (9,-9)**@{-}; (-9,-9)**@{-};
(-9,3); (-9,5)**@{-}; (-3,5)**@{-}; (-3,3)**@{-}; (-9,3)**@{-};
(-9,7); (-9,9)**@{-}; (-7,9)**@{-}; (-7,7)**@{-}; (-9,7)**@{-};
(-5,7); (-5,9)**@{-}; (-3,9)**@{-}; (-3,7)**@{-}; (-5,7)**@{-};
(9,3); (9,9)**@{-}; (3,9)**@{-}; (3,3)**@{-}; (9,3)**@{-};
\endxy
\]

Each numbered pattern corresponds to a numbered partition of
\(C\times C\) in the obvious way.  To make sure that we agree on
this, we give the above example a numbering.

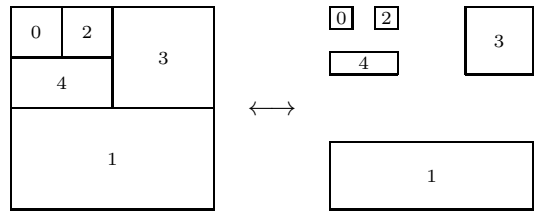
\[
\xy (-13.5,-13.5); (-13.5,13.5)**@{-}; (13.5,13.5)**@{-}; (13.5,-13.5)**@{-};
(-13.5,-13.5)**@{-};
(-13.5,0); (13.5,0)**@{-}; (0,0); (0,13.5)**@{-};
(-13.5,6.75); (0,6.75)**@{-};
(-6.75,6.75); (-6.75,13.5)**@{-};
(0,-6.75)*{\scriptstyle1};
(6.75,6.75)*{\scriptstyle3};
(-6.75,3.375)*{\scriptstyle4};
(-10.125,10,125)*{\scriptstyle0};
(-3.375,10,125)*{\scriptstyle2};
\endxy
\quad
\longleftrightarrow
\quad
\xy
(-13.5,-13.5); (-13.5,-4.5)**@{-}; (13.5,-4.5)**@{-}; 
(13.5,-13.5)**@{-}; (-13.5,-13.5)**@{-};
(0,-9)*{\scriptstyle1};
(-13.5,4.5); (-13.5,7.5)**@{-}; (-4.5,7.5)**@{-}; 
(-4.5,4.5)**@{-}; (-13.5,4.5)**@{-};
(-9,6)*{\scriptstyle4};
(-13.5,10.5); (-13.5,13.5)**@{-}; (-10.5,13.5)**@{-}; 
(-10.5,10.5)**@{-}; (-13.5,10.5)**@{-};
(-12,12)*{\scriptstyle0};
(-7.5,10.5); (-7.5,13.5)**@{-}; (-4.5,13.5)**@{-}; 
(-4.5,10.5)**@{-}; (-7.5,10.5)**@{-};
(-6,12)*{\scriptstyle2};
(13.5,4.5); (13.5,13.5)**@{-}; (4.5,13.5)**@{-}; 
(4.5,4.5)**@{-}; (13.5,4.5)**@{-};
(9,9)*{\scriptstyle3};
\endxy
\]

\subsection{Homeomorphisms from numbered pattern
pairs}\mylabel{HomeoFromPatt}

Let \(P_1\) and \(P_2\) be numbered patterns with the same number
\(n\) of rectangles in each.  Let the rectangles in \(P_1\) be
\(R_0\), \(R_1\), \dots, \(R_{n-1}\) where the subscripts reflect
the numbering, and let the rectangles in \(P_2\) be \(T_0\),
\(T_1\), \dots, \(T_{n-1}\) with the same comment.  We will get a
self homeomorphism \(h=h(P_1, P_2)\) of \(C\times C\) from the pair
\((P_1, P_2)\).

We let \(A(R_i)\) and \(A(T_i)\) be the closed and open sets in
\(C\times C\) associated, respectively, to \(R_i\) and \(T_i\), and
we let \(h\) take \(A(R_i)\) onto \(A(T_i)\) affinely so as to
preserve the orientation in each coordinate.  By this we mean that
the restriction of \(h\) to \(A(R_i)\) is the restriction of the
unique self homeomorphism of \(\R^2\) of the form \((x,y)\mapsto
(a+3^jx, b+3^ky)\) that maps \(A(R_i)\) onto \(A(T_i)\).  Thus the
four corners of \(A(R_i)\), the lower left, lower right, upper
right, upper left, are each carried to the corner of \(A(T_i)\) with
the same description.  Doing this for each \(i\) in \(\{0,1,\ldots,
n-1\}\) defines a homeomorphism \(h(P_1, P_2)\) from \(C\times C\)
to itself.

We picture such homeomorphisms by giving a pair of numbered patterns
with an arrow from the numbered pattern describing the domain to the
numbered pattern describing the range.  The pair below represents an
element called the ``baker's map.''  This particular map will be
discussed later in this paper.

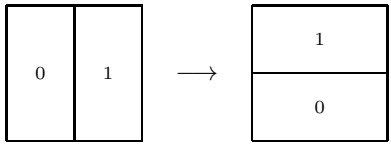
\[
\xy
(-9,-9); (-9,9)**@{-}; (9,9)**@{-}; (9,-9)**@{-}; (-9,-9)**@{-};
(0,-9); (0,9)**@{-};
(-4.5,0)*{\scriptstyle0};
(4.5,0)*{\scriptstyle1};
\endxy
\quad
\longrightarrow
\quad
\xy
(-9,-9); (-9,9)**@{-}; (9,9)**@{-}; (9,-9)**@{-}; (-9,-9)**@{-};
(-9,0); (9,0)**@{-};
(0,-4.5)*{\scriptstyle0};
(0,4.5)*{\scriptstyle1};
\endxy
\]

\subsection{The group \protect\(2V\protect\)}

The group \(2V\) is the set of all self homeomorphisms of \(C\times
C\) of the form \(h(P_1, P_2)\) where \(P_1\) and \(P_2\) are
numbered patterns with the same number of rectangles.  The group
operation is composition.  Closure under inversion is immediate and
closure under composition is a pleasant exercise.  Those that wish
to defer the exercise can wait until Section \ref{WordPairMult}
where techniques for multiplying elements will be given.

The following is clear.

\begin{lemma}\mylabel{BigInfinite} The group \(2V\) is countably
infinite and contains all finite groups.  \end{lemma}

\section{Simplicity of the commutator
subgroup}\mylabel{SimplCommSect}

It is usually the case that a sufficiently transitive permutation
group generated by elements ``of small support'' has a simple
commutator subgroup.  See \cite{higman:simpl} or
\cite{epstein:simpl}.  This section contains the arguments to show
the simplicity of the commutator subgroup of \(2V\).  This section
is self contained since it is often just as easy to apply the
techniques of \cite{higman:simpl} and \cite{epstein:simpl} as it is
to quote them.

Later in the paper we will know that \(2V\) equals its commutator
subgroup.  Thus we will have proven that \(2V\) is simple.

\subsection{Transitivity}

We establish that \(2V\) is sufficiently ``transitive.''

\begin{lemma}\mylabel{SuffTrans} Let \(K\) be a closed, proper
subset of \(C\times C\) and let \(U\) be a non-empty open set in
\(C\times C\).  Then there is an element \(h\) of \(2V\) so that
\(h(K)\) is contained in \(U\).  \end{lemma}

\begin{proof} Since \(K\) is closed and not all of \(C\times C\),
there is a pattern \(P_1\) where not all the rectangles are needed
to get a set \(S\) of rectangles whose corresponding closed and open
sets in \(C\times C\) cover \(K\).  There is a pattern \(P_2\) with
more than one rectangle and with a rectangle \(R\) so that \(A(R)\)
is contained in \(U\).  If \(n\) is the number of rectangles in
\(S\), then we can apply \(n-1\) subdivision operations at random to
the rectangle \(R\) to create a new pattern \(P'_2\) containing
\(n\) different rectangles that are contained in \(R\) whose
corresponding sets in \(C\) are thus contained in \(U\).  Using
further subdivisions if necessary in some rectangle not in \(S\) or
in \(R\), we can get the number of rectangles in \(P_1\) and
\(P'_2\) to be the same.  Now an element of \(2V\) can be built that
carries the rectangles in \(S\) into \(R\).  \end{proof}

\subsection{Small generators}

Our next step is to argue that \(2V\) is generated by elements ``of
small support.''

The word small is misleading.  All that we will need from the notion
of small is that a composition of a fixed number of elements of
sufficiently small support yields an element whose support is a
closed, proper subset of \(C\times C\).  Our notion of small will
also cooperate equally well with \(C^n\) and \(C^\N\).

We will measure size in one coordinate only.  Since \(C\) has
measure 0, we will not measure in \(C\) itself, but in \(I\).  We
use the notion of patterns in \(I\) in a manner parallel to patterns
in \(I^2\).  The unit interval \(I=[0,1]\) is the trivial pattern
and all other patterns are derived from the trivial pattern by
dividing intervals at their midpoints.

Given an \(\epsilon > 0\), we say that a set in \(C\) is \emph{of
size no more than} \(\epsilon\) if it is in a collection of disjoint
closed and open sets in \(C\) corresponding to a set of intervals in
a pattern on \(I\) whose lengths sum to no more than \(\epsilon\).
Now we say that a set in \(C\times C\) is \emph{of size no more
than} \(\epsilon\) if the projection of the set to the first
coordinate is of size no more than \(\epsilon\).  Given
\(\epsilon>0\), let \(B_\epsilon\) be the set of those elements of
\(2V\) that have support of size no more than \(\epsilon\).

\begin{prop}\mylabel{SmallSupsGen} Given \(\epsilon>0\), the group
\(2V\) is generated by \(B_\epsilon\).  \end{prop}

\begin{proof} Let \(A_{0}\) be the elements in \(2V\) with support
in the left half of \(C\times C\), and let \(A_{1}\) be the elements
with support in the right half of \(C\times C\).  The elements of
\(A_{0}\) and \(A_1\) have support of size no more than \(\frac12\).

Let \(h\) be in \(2V\).  We will show first that \(h\) is the
composition of elements of \(B_\epsilon\cup A_0\cup A_1\).  The
claimed result will follow from this because an identical exercise
will then show that elements in each \(A_i\) will be compositions of
elements in \(B_\epsilon\cup A_{i0}\cup A_{i1}\) where \(A_{i0}\)
and \(A_{i1}\) are the sets of elements whose supports are contained
in, respectively, the left and right halves of the half of \(C\times
C\) containing the support of the elements of \(A_i\).  Thus
elements of the \(A_{ij}\) have supports of size no more than
\(\frac14\).  Continuing in this way, we get \(h\) to be a
composition of elements in \(B_\epsilon\).

We will alter \(h\) and then work with the alterations.  To simplify
the notation, we will use \(h\) to refer to the current alteration.
The alterations will be to compose \(h\) with elements of \(A_0\)
and \(A_1\), and the goal will be to reduce \(h\) to an element of
\(B_\epsilon\).

We may assume that \(h\) is not trivial and thus \(x\ne h(x)\) for
some \(x\).  Let us assume for now that \(x\) is in the left half of
\(C\times C\).  There are two disjoint closed and open sets \(E\)
and \(F\) corresponding to rectangles in patterns with \(x\in E\),
\(h(x)\in F\) so that \(h\) carries \(E\) affinely onto \(F\).  We
can choose \(E\) to be contained in the left half of \(C\times C\).
We do not care about the location of \(F\).

At the expense of making the sets smaller, we can find a pattern so
that both \(E\) and \(F\) correspond to rectangles in the pattern.
This allows us to build an element \(h_1\) that agrees with \(h\) on
\(E\), and that does nothing but interchange \(E\) and \(F\) so that
\(\supp(h_1)\) is in \(E\cup F\).  Using the uniform continuity of
\(h\), we can choose \(E\) small enough so that \(h_1\) is in
\(B_\epsilon\).  Because \(E\) and \(F\) are disjoint, we note that
\(h^{-1}_1h\) fixes the non-empty open set \(E\) and that
\(\supp(h^{-1}_1h)\) is contained in \(\supp(h)\).  We replace \(h\)
by \(h^{-1}_1h\).

Using Lemma \ref{SuffTrans}, we can build an element \(g\) of
\(A_{0}\) that carries the subset \(E\) onto the complement in the
left half of \(C\times C\) of an open set \(U\) whose closure has
size no more than \(\frac\epsilon2\).  Now
\(\supp(h^g)=g(\supp(h))\), so we know that the size of
\(\supp(h^g)\) in the left half of \(C\times C\) is of size no more
than \(\frac\epsilon2\).  We replace \(h\) by \(h^g\).

We now move our attention to the right half of \(C\times C\).  If
there is some \(x\) there with \(x\ne h(x)\), then we repeat the
above operations except that the conjugating element \(g\) is from
\(A_{1}\).  The result is an \(h\) whose support in the right half
of \(C\times C\) is of size no more than \(\frac\epsilon2\).  Of the
two steps in the modification, the first (replacing \(h\) by
\(h^{-1}_1h\)) does not increase support, and the second
(conjugating by \(g\)) only modifies the support in the right half
of \(C\times C\).  Thus these modifications will not affect the
fact that the support of the new \(h\) in the left half of \(C\times
C\) is of size no more than \(\frac\epsilon2\).  This reduces the
support of \(h\) to have size no more than \(\epsilon\).
\end{proof}

\subsection{Simplicity of the commutator subgroup} The ideas in the
next lemma are to be found in \cite{higman:simpl} and
\cite{epstein:simpl}.

\begin{prop}\mylabel{SimplComm} The commutator subgroup of \(2V\) is
simple.  \end{prop}

\begin{proof} Let \(W\) be the commutator subgroup of \(2V\).  Let
\(j\) be a non-trivial element of \(W\) and let \(N\) be the normal
closure of \(j\) in \(W\).  We must show that elements of \(2V\)
commute modulo \(N\).  It suffices to show that any two generators
of \(2V\) commute modulo \(N\).  Thus we let \(h\) and \(g\) be
elements of \(2V\) with very small supports.

There is some open and closed set \(E\) in the support of \(j\) so
that \(E\) and \(j(E)\) are disjoint.  If some element \(k\) of
\(W\) carries the union of the supports of \(g\) and \(h\) into
\(E\), then \[\supp\left(h^{(j^{k^{-1}})}\right) =k^{-1}jk
(\supp(h)) \subseteq k^{-1}j(E) \] is disjoint from
\(\supp(g)\subseteq k^{-1}(E)\), and \(g\) and \(h\) commute modulo
\(N\).  Thus we are done if we find such an element \(k\).

Let \(X\) be the union of the supports of \(h\) and \(g\).  By Lemma
\ref{SuffTrans}, there is an element \(k'\) in \(2V\) that carries
\(X\) into \(E\).  Since this can be done by interchanging
rectangles, we can keep the support of \(k'\) contained in a small
neighborhood of \(X\cup E\).  Since this can be kept very small,
there is an \(r\) in \(2V\) that carries \(X\cup E\) to a set
disjoint from \(X\cup E\).  Thus it is seen that the element
\(k=[k',r] = k'\left((k')^{-1}\right)^r\) of \(W\) agrees with
\(k'\) on \(X\).  This completes the proof.  \end{proof}

\section{Groups related to \protect\(2V\protect\)}

\subsection{The groups \protect\(nV\protect\) and
\protect\(\protect\omega V\protect\)}

For each \(n\) there is a group \(nV\) that acts on the product
\(C^n\).  There is also a group \(\omega V\) that acts on the
countably infinite product \(C^\N\).  It is the ascending union of
the \(nV\).  Patterns for \(nV\) would consist of rectangular solids
in \(I^n\).  These are easiest to code by \(n\)-tuples of finite
(and possibly empty) words in the alphabet \(\{0,1\}\).  The
\(n\)-tuple of empty words corresponds to \(I^n\).  If a rectangular
solid is given by the \(n\)-tuple \((w_0, w_1, \ldots, w_{n-1})\),
then a division into two solids by a cut perpendicular to coordinate
axis \(i\) would result in the two solids given by replacing \(w_i\)
by \(w_i0\) or \(w_i1\).

For \(\omega V\), solids in \(I^\N\) would be given by infinite
sequences of finite words in \(\{0,1\}\) such that all but finitely
many of the words in the sequence are the empty word.  A division of
a solid would be similar to a division in \(I^n\).

The group \(1V\) is the Thompson group \(V\) as described in
\cite{CFP}.  We show below that \(V\) and \(2V\) are not isomorphic.
The groups \(2V\) and \(\omega V\) are not isomorphic since one is
finitely generated and the other is not.  It is not known whether
any other pairs of the \(nV\) with \(n>1\) or \(\omega V\) are
isomorphic to each other, but it would be pleasant if the answer
were no.

It is clear that some of the results of this paper and
\cite{brin:hd4} should apply to the \(nV\) and \(\omega V\), but it
is not clear how many.  Lemma \ref{BigInfinite} clearly applies.
The simplicity results of this paper do not depend on finding a
presentation, but they do depend on finding enough relations to
prove that the abelianization is trivial.  At this point, one can
only say that it is believable that this can be done for the \(nV\)
and \(\omega V\).  It can be said that the results of Section
\ref{SimplCommSect} apply and that \(nV\) and \(\omega V\) possess
simple commutator subgroups.  Whether a calculation of a full,
finite presentation can be carried out for all these groups comes
under other descriptions.  For the \(nV\) with \(n\) finite, it is
to be hoped for.  The group \(\omega V\) is another matter since it
is not finitely generated.  For \(\omega V\) the hope would be to
write down any reasonable presentation.

\subsection{Other related groups} In \cite{higman:fpsimpl}, groups
are defined that are called \(G_{n,r}\) there, but which we will
relable as \(V_{k,r}\) here to conform to the letter \(V\) used in
\cite{CFP} and to avoid conflict with our use of \(n\).  Groups
\(nV_{k,r}\) can be defined to parallel the \(V_{k,r}\) so that
\(1V_{k,r}=V_{k,r}\).  The group \(nV_{k,r}\) would consist of self
homeomorphisms of \(r\) copies of \(C^n\).  To obtain the defining
patterns of these homeomorphisms, the basic inductive step would be
to divide a rectangular solid into \(k\) congruent pieces by cuts
perpendicular to a particular axis.  We will not treat the
\(nV_{k.r}\) here.  It should be noted, that not all the \(V_{k,r}\)
are simple, but in \cite{higman:fpsimpl} it is shown they have
infinte, simple, finitely presented commutator subgroups that are of
index 1 when \(k\) is even and index 2 when \(k>2\) is odd.

\subsection{The containing group
\protect\(\protect\widehat{2V}\protect\) of \protect\(2V\protect\)}

The group \(2V\) is contained in a ``larger'' finitely presented
group \(\widehat{2V}\).  The relation between \(2V\) and
\(\widehat{2V}\) is similar to the relation between \(BV\) and
\(\widehat{BV}\) in \cite{brin:bv}.  In \cite{brin:bv}, the group
\(\widehat{BV}\) is analyzed first and this analysis is then used to
analyze its subgroup \(BV\).  This outline will be partly followed
here.  We briefly discuss why and why only partly.

Once generators are found for \(2V\), the word problem is
theoretically solvable.  The elements of \(2V\) are specific
permutations of \(C\times C\) and a word in specific elements is
trivial if and only if the resulting permutation is trivial.
However, there is a large difference between theory and practice and
we need a more practical way to identify the trivial element.

The group \(BV\) involves braids, and identifying the trivial
element in \(BV\) is even less trivial.  Triviality in \(BV\) is
detected by its setting as a subgroup of \(\widehat{BV}\) which has
a known, finite presentation and a strong normal form (Lemma 10.3 of
\cite{brin:bv}).  While \(\widehat{2V}\) has a known, finite
presentation, it does not have as nice a normal form.  This is
discussed in \cite{brin:hd4}.  Thus we need another technique for
identifying trivial elements.  This is solved in \cite{brin:hd4} by
finding a less nice normal form for elements of \(\widehat{2V}\) and
an algorithm for achieving the form.

Even before we need to understand trivial elements, we need to be
able to multiply easily in \(2V\).  It turns out that
\(\widehat{2V}\) has a very easy multiplication since it is a group
of fractions of a very well behaved monoid \(\Pi\) of positive
elements with an even easier multiplication.  Applying the
multiplication in \(\Pi\) to \(2V\) will be discussed more fully
below.  (See Section \ref{WordPairMult}.)

One last reason for discussing \(\widehat{2V}\) is that the
advertised resemblance to the little cubes operads of
\cite{may:geom+loop} and \cite{boardman+vogt} is seen most strongly
in the monoid \(\Pi\).

\subsection{The positive monoid of
\protect\(\protect\widehat{2V}\protect\)}  

We start with the monoid \(\Pi\).  Elements of the monoid will
correspond to certain ``numbered sequences of patterns.''  We make
this specific by starting with an infinite sequence of pairwise
disjoint unit squares \((S_0, S_1, \ldots, )\) in the plane.  To be
very specific, we can take them in the upper half plane so that
\(S_i\) intersects the \(x\)-axis in the interval \([2i, 2i+1]\).
Now a sequence of patterns (we will number them shortly) is an
infinite sequence \((P_0, P_1, \ldots)\) of patterns where \(P_i\)
is thought of as a pattern in \(S_i\) and only finitely many of the
\(P_i\) are not trivial.

A numbering of such a sequence of patterns is a one-to-one
correspondence between \(\N\) and the rectangles in the sequence for
which there are \(j\) and \(k\) in \(\N\) so that \(i>k\) implies
that the pattern \(P_i\) that is applied to \(S_i\) is the trivial
pattern and that the number of the rectangle consisting of \(S_i\)
is \(i+j\).  In words, the numbering eventually becomes a
consecutive numbering, from left to right, of unsubdivided unit
squares.  An example is given below where it is assumed that all
squares not pictured are numbered consecutively from left to right
starting with 10.  Thus in this example, \(k=3\) and \(j=5\) work
for the required restriction.

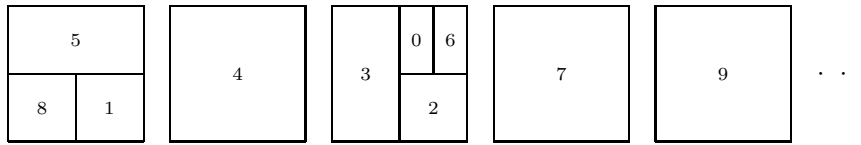
\[
\xy
(-9,-9); (9,-9)**@{-}; (9,9)**@{-}; (-9,9)**@{-}; (-9,-9)**@{-};
(-9,0); (9,0)**@{-};
(0,-9); (0,0)**@{-};
(0,4.5)*{\scriptstyle5};
(-4.5,-4.5)*{\scriptstyle8};
(4.5,-4.5)*{\scriptstyle1};
\endxy
\quad
\xy
(-9,-9); (9,-9)**@{-}; (9,9)**@{-}; (-9,9)**@{-}; (-9,-9)**@{-};
(0,0)*{\scriptstyle4};
\endxy
\quad
\xy
(-9,-9); (9,-9)**@{-}; (9,9)**@{-}; (-9,9)**@{-}; (-9,-9)**@{-};
(0,-9); (0,9)**@{-};
(0,0); (9,0)**@{-};
(4.5,0); (4.5,9)**@{-};
(-4.5,0)*{\scriptstyle3};
(4.5,-4.5)*{\scriptstyle2};
(2.25,4.5)*{\scriptstyle0};
(6.75,4.5)*{\scriptstyle6};
\endxy
\quad
\xy
(-9,-9); (9,-9)**@{-}; (9,9)**@{-}; (-9,9)**@{-}; (-9,-9)**@{-};
(0,0)*{\scriptstyle7};
\endxy
\quad
\xy
(-9,-9); (9,-9)**@{-}; (9,9)**@{-}; (-9,9)**@{-}; (-9,-9)**@{-};
(0,0)*{\scriptstyle9};
\endxy
\quad
\xy
(0,0)*{\cdot};
(3,0)*{\cdot};
(6,0)*{\cdot};
\endxy
\]

If we let \(X\) be the union of the \(S_i\), then a numbered
sequence of patterns \((P_0, P_1, \ldots)\) describes a continuous
function from \(X\) to itself.  The function simply takes the square
\(S_i\), affinely and with orientation of each coordinate preserved,
to rectangle \(i\) in the numbered sequence of patterns.  Under
composition of functions, this is a submonoid of the monoid
\(C^0(X)\) of continuous functions from \(X\) to itself.  The
identity element is given by the trivial sequence in which every
\(S_i\) is a rectangle and the number of \(S_i\) is \(i\).

If \(h\) is an element of the monoid described by a sequence of
patterns \(P=(P_0, P_1, \ldots)\), then we can think of \(h\) as
carrying the trivial sequence to the sequence \(P\).  Thus we can be
less efficient and think of two pieces of information as describing
\(h\): the sequence \(P\) and the trivial sequence.  The trivial
sequence describes a structure of the domain of \(h\) and the
sequence \(P\) describes the structure of the range of \(h\).  This
lack of efficiency anticipates what we need to describe elements of
the full group \(\widehat{2V}\).

It is easy to argue that this monoid is cancellative and has common
right multiples (if composition proceeds from right to left) and
thus has a group of fractions by Ore's theorem (see Theorem 1.23 of
\cite{cliff+prest:I}).  However, this is not necessary.  We can
``invert'' the elements by representing them as homeomorphisms
rather than continuous functions.  This version with invertible
elements will be what we call the monoid \(\Pi\).

\subsection{The group \protect\(\protect\widehat{2V}\protect\)}

Each square \(S_i\) contains a copy of \(C\times C\) in the obvious
way, and we let \(Y\) be the union of all these copies of \(C\times
C\).  The numbered sequence of patterns gives a numbering of closed
and open sets in \(Y\) that correspond to the rectangles in the
sequence of pattens.  Now a homeomorphism from \(Y\) to \(Y\) comes
out of a numbered sequence by taking the entire copy of \(C\times
C\) in \(S_i\) to the closed and open set numbered \(i\) by an
affine map that preserves orientation in each coordinate.  This is
done exactly as described in Section \ref{HomeoFromPatt}.  

We will use \(\Pi\) to denote this monoid of self homeomorphisms of
\(Y\).  The group of fractions that we want is the group generated
by the elements of \(\Pi\) and their inverses.  That this
corresponds to the group of fractions construction follows from
Problem 3 on Page 37 and Theorem 1.24 of \cite{cliff+prest:I}.

The group \(\widehat{2V}\) is the group just described.  Elements
can be represented by pairs of numbered sequences of patterns.  If
\(P\) and \(Q\) are patterns representing elements \(g\) and \(h\),
respectively, in \(\Pi\), then each is a homeomorphism from \(Y\) to
\(Y\).  Each has its domain structure described by the trivial
sequence of patterns.  The range structure for \(g\) is described by
\(P\) and the range structure for \(h\) is described by \(Q\).  If
we consider the element \(gh^{-1}\) of \(\widehat{2V}\), then the
resulting homeomorphism from \(Y\) to itself is described by using
\(Q\) to determine the structure of the domain and by using \(P\) to
determine the structure of the range.  For each \(i\), the
homeomorphism \(gh^{-1}\) takes the set corresponding to rectangle
\(i\) in \(Q\) to the set corresponding to the rectangle \(i\) in
\(P\).

If we use the pair \((P,Q)\) to represent the element \(gh^{-1}\) in
\(\widehat{2V}\), then we can think of this as a map from \(Q\) to
\(P\).  We use the right element of \((P,Q)\) to represent the
pattern for the domain so that right to left composition reads
nicely as in \((P,Q)(Q,R)=(P,R)\).

One can also define groups \(\widehat{nV}\) and \(\widehat{\omega
V}\) in a similar manner, but we have nothing to say about
these groups.

In the next few paragraphs, we will discuss the product in \(\Pi\),
we will describe some elements of \(\Pi\) that clearly generate
\(\Pi\), and we will describe some relations satisfied by those
elements.  The relations that we give suffice to present \(\Pi\) but
this will be less clear and a proof of this fact will be given in
\cite{brin:hd4}.

\subsection{Multiplication in the positive monoid of
\protect\(\protect\widehat{2V}\protect\)}\mylabel{HatTwoVMonMult}

If a sequence of patterns \(P=(P_0, P_1, \ldots)\) represents an
element \(h\) in the positive monoid \(\Pi\) of \(\widehat{2V}\) and
\(Q=(Q_0, Q_1, \ldots)\) represents \(g\in \Pi\), then the sequence
of patterns \(PQ\) for \(hg\) (recall that \(g\) is applied first)
is gotten by pasting the pattern \(Q_i\) affinely into the rectangle
\(i\) of \(P\) for each \(i\in\N\).  The numbering of the resulting
pattern is that of \(Q\).  That is, when the pattern \(Q_i\) is
pasted into rectangle \(i\) of \(P\), the rectangle numbers of \(Q\)
are pasted in along with the pattern.  We illustrate this below
where the squares not pictured in \(P\) are not subdivided and as
rectangles they are numbered consecutively, left to right, from 7.
Similarly in \(Q\) the corresponding number is 14, and in \(PQ\) the
corresponding number is 17.

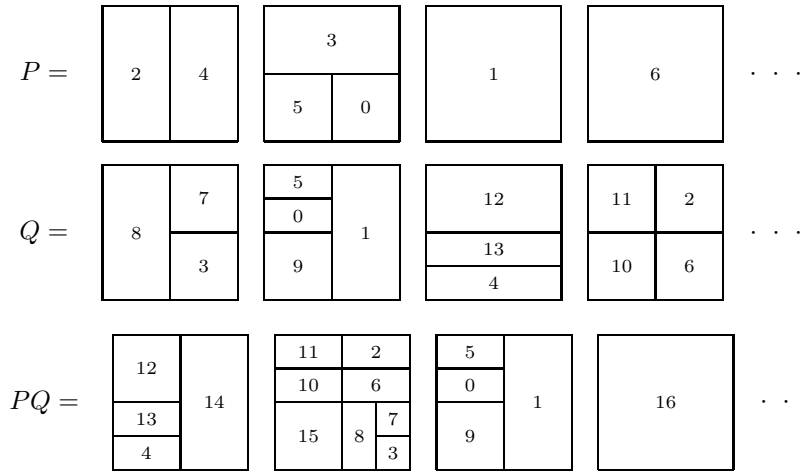
\[
P=\quad
\xy
(-9,-9); (-9,9)**@{-}; (9,9)**@{-}; (9,-9)**@{-}; (-9,-9)**@{-};
(0,-9); (0,9)**@{-};
(-4.5,0)*{\scriptstyle2};
(4.5,0)*{\scriptstyle4};
\endxy
\quad
\xy
(-9,-9); (-9,9)**@{-}; (9,9)**@{-}; (9,-9)**@{-}; (-9,-9)**@{-};
(-9,0); (9,0)**@{-};
(0,0); (0,-9)**@{-};
(0,4.5)*{\scriptstyle3};
(-4.5,-4.5)*{\scriptstyle5};
(4.5,-4.5)*{\scriptstyle0};
\endxy
\quad
\xy
(-9,-9); (-9,9)**@{-}; (9,9)**@{-}; (9,-9)**@{-}; (-9,-9)**@{-};
(0,0)*{\scriptstyle1};
\endxy
\quad
\xy
(-9,-9); (-9,9)**@{-}; (9,9)**@{-}; (9,-9)**@{-}; (-9,-9)**@{-};
(0,0)*{\scriptstyle6};
\endxy
\quad
\xy
(-3,0)*{\cdot};
(0,0)*{\cdot};
(3,0)*{\cdot};
\endxy
\]

\[
Q=\quad
\xy
(-9,-9); (-9,9)**@{-}; (9,9)**@{-}; (9,-9)**@{-}; (-9,-9)**@{-};
(0,-9); (0,9)**@{-};
(0,0); (9,0)**@{-};
(-4.5,0)*{\scriptstyle8};
(4.5,4.5)*{\scriptstyle7};
(4.5,-4.5)*{\scriptstyle3};
\endxy
\quad
\xy
(-9,-9); (-9,9)**@{-}; (9,9)**@{-}; (9,-9)**@{-}; (-9,-9)**@{-};
(0,-9); (0,9)**@{-};
(0,0); (-9,0)**@{-};
(-9,4.5); (0,4.5)**@{-};
(4.5,0)*{\scriptstyle1};
(-4.5,2.25)*{\scriptstyle0};
(-4.5,6.75)*{\scriptstyle5};
(-4.5,-4.5)*{\scriptstyle9};
\endxy
\quad
\xy
(-9,-9); (-9,9)**@{-}; (9,9)**@{-}; (9,-9)**@{-}; (-9,-9)**@{-};
(-9,0); (9,0)**@{-};
(-9,-4.5); (9,-4.5)**@{-};
(0,-2.25)*{\scriptstyle13};
(0,-6.75)*{\scriptstyle4};
(0,4.5)*{\scriptstyle12};
\endxy
\quad
\xy
(-9,-9); (-9,9)**@{-}; (9,9)**@{-}; (9,-9)**@{-}; (-9,-9)**@{-};
(0,-9); (0,9)**@{-};
(-9,0); (9,0)**@{-};
(-4.5,4.5)*{\scriptstyle11};
(-4.5,-4.5)*{\scriptstyle10};
(4.5,4.5)*{\scriptstyle2};
(4.5,-4.5)*{\scriptstyle6};
\endxy
\quad
\xy
(-3,0)*{\cdot};
(0,0)*{\cdot};
(3,0)*{\cdot};
\endxy
\]

\[
PQ=\quad
\xy
(-9,-9); (-9,9)**@{-}; (9,9)**@{-}; (9,-9)**@{-}; (-9,-9)**@{-};
(0,-9); (0,9)**@{-};
(-9,0); (0,0)**@{-};
(-9,-4.5); (0,-4.5)**@{-};
(-4.5,4.5)*{\scriptstyle12};
(-4.5,-2.25)*{\scriptstyle13};
(-4.5,-6.75)*{\scriptstyle4};
(4.5,0)*{\scriptstyle14};
\endxy
\quad
\xy
(-9,-9); (-9,9)**@{-}; (9,9)**@{-}; (9,-9)**@{-}; (-9,-9)**@{-};
(-9,0,); (9,0)**@{-};
(0,0); (0,-9)**@{-};
(4.5,-9); (4.5,0)**@{-};
(4.5,-4.5); (9,-4.5)**@{-};
(-9,4.5); (9,4.5)**@{-};
(0,0); (0,9)**@{-};
(-4.5,-4.5)*{\scriptstyle15};
(2.25,-4.5)*{\scriptstyle8};
(6.75,-2.25)*{\scriptstyle7};
(6.75,-6.75)*{\scriptstyle3};
(-4.5,2.25)*{\scriptstyle10};
(-4.5,6.75)*{\scriptstyle11};
(4.5,2.25)*{\scriptstyle6};
(4.5,6.75)*{\scriptstyle2};
\endxy
\quad
\xy
(-9,-9); (-9,9)**@{-}; (9,9)**@{-}; (9,-9)**@{-}; (-9,-9)**@{-};
(0,-9); (0,9)**@{-};
(0,0); (-9,0)**@{-};
(-9,4.5); (0,4.5)**@{-};
(4.5,0)*{\scriptstyle1};
(-4.5,2.25)*{\scriptstyle0};
(-4.5,6.75)*{\scriptstyle5};
(-4.5,-4.5)*{\scriptstyle9};
\endxy
\quad
\xy
(-9,-9); (-9,9)**@{-}; (9,9)**@{-}; (9,-9)**@{-}; (-9,-9)**@{-};
(0,0)*{\scriptstyle16};
\endxy
\quad
\xy
(-3,0)*{\cdot};
(0,0)*{\cdot};
(3,0)*{\cdot};
\endxy
\]

\subsection{Elements of the positive monoid of
\protect\(\protect\widehat{2V}\protect\)}\mylabel{MonoidGens}

For \(i\ge0\), let \(v_i\) be as pictured below.

\[
v_i=
\xy
(-6,-6); (-6,6)**@{-}; (6,6)**@{-}; (6,-6)**@{-}; (-6,-6)**@{-};
(0,0)*{\scriptstyle0};
\endxy
\quad
\xy
(-6,-6); (-6,6)**@{-}; (6,6)**@{-}; (6,-6)**@{-}; (-6,-6)**@{-};
(0,0)*{\scriptstyle1};
\endxy
\quad
\xy
(-3,0)*{\cdot};
(0,0)*{\cdot};
(3,0)*{\cdot};
\endxy
\quad
\xy
(-6,-6); (-6,6)**@{-}; (6,6)**@{-}; (6,-6)**@{-}; (-6,-6)**@{-};
(0,0)*{\scriptstyle{i-1}};
\endxy
\quad
\xy
(-6,-6); (-6,6)**@{-}; (6,6)**@{-}; (6,-6)**@{-}; (-6,-6)**@{-};
(0,-6);(0,6)**@{-};
(-3,0)*{\scriptstyle{i}};
(3,0)*{\scriptstyle{i+1}};
\endxy
\quad
\xy
(-6,-6); (-6,6)**@{-}; (6,6)**@{-}; (6,-6)**@{-}; (-6,-6)**@{-};
(0,0)*{\scriptstyle{i+2}};
\endxy
\quad
\xy
(-3,0)*{\cdot};
(0,0)*{\cdot};
(3,0)*{\cdot};
\endxy
\]

In the above picture, each square \(S_j\) with \(j\ne i\) has the
trivial pattern, each square \(S_j\) with \(j<i\) is numbered \(j\)
and each square \(S_j\) with \(j>i\) is numbered \(j+1\).

For \(i\ge 0\), let \(h_i\) be as pictured below.

\[
h_i=
\xy
(-6,-6); (-6,6)**@{-}; (6,6)**@{-}; (6,-6)**@{-}; (-6,-6)**@{-};
(0,0)*{\scriptstyle0};
\endxy
\quad
\xy
(-6,-6); (-6,6)**@{-}; (6,6)**@{-}; (6,-6)**@{-}; (-6,-6)**@{-};
(0,0)*{\scriptstyle1};
\endxy
\quad
\xy
(-3,0)*{\cdot};
(0,0)*{\cdot};
(3,0)*{\cdot};
\endxy
\quad
\xy
(-6,-6); (-6,6)**@{-}; (6,6)**@{-}; (6,-6)**@{-}; (-6,-6)**@{-};
(0,0)*{\scriptstyle{i-1}};
\endxy
\quad
\xy
(-6,-6); (-6,6)**@{-}; (6,6)**@{-}; (6,-6)**@{-}; (-6,-6)**@{-};
(-6,0);(6,0)**@{-};
(0,-3)*{\scriptstyle{i}};
(0,3)*{\scriptstyle{i+1}};
\endxy
\quad
\xy
(-6,-6); (-6,6)**@{-}; (6,6)**@{-}; (6,-6)**@{-}; (-6,-6)**@{-};
(0,0)*{\scriptstyle{i+2}};
\endxy
\quad
\xy
(-3,0)*{\cdot};
(0,0)*{\cdot};
(3,0)*{\cdot};
\endxy
\]

In the above picture, each square \(S_j\) with \(j\ne i\) has the
trivial pattern, each square \(S_j\) with \(j<i\) is numbered \(j\)
and each square \(S_j\) with \(j>i\) is numbered \(j+1\).

For \(i\ge 0\), let \(\sigma_i\) be as pictured below.

\[
\sigma_i=
\xy
(-6,-6); (-6,6)**@{-}; (6,6)**@{-}; (6,-6)**@{-}; (-6,-6)**@{-};
(0,0)*{\scriptstyle0};
\endxy
\quad
\xy
(-6,-6); (-6,6)**@{-}; (6,6)**@{-}; (6,-6)**@{-}; (-6,-6)**@{-};
(0,0)*{\scriptstyle1};
\endxy
\quad
\xy
(-3,0)*{\cdot};
(0,0)*{\cdot};
(3,0)*{\cdot};
\endxy
\quad
\xy
(-6,-6); (-6,6)**@{-}; (6,6)**@{-}; (6,-6)**@{-}; (-6,-6)**@{-};
(0,0)*{\scriptstyle{i-1}};
\endxy
\quad
\xy
(-6,-6); (-6,6)**@{-}; (6,6)**@{-}; (6,-6)**@{-}; (-6,-6)**@{-};
(0,0)*{\scriptstyle{i+1}};
\endxy
\quad
\xy
(-6,-6); (-6,6)**@{-}; (6,6)**@{-}; (6,-6)**@{-}; (-6,-6)**@{-};
(0,0)*{\scriptstyle{i}};
\endxy
\quad
\xy
(-6,-6); (-6,6)**@{-}; (6,6)**@{-}; (6,-6)**@{-}; (-6,-6)**@{-};
(0,0)*{\scriptstyle{i+2}};
\endxy
\quad
\xy
(-3,0)*{\cdot};
(0,0)*{\cdot};
(3,0)*{\cdot};
\endxy
\]

In the above picture, every \(S_j\) has the trivial pattern and
every \(S_j\) with \(j\notin\{i,i+1\}\) is numbered \(j\).

We make some remarks about multiplying by these elements on the
right.  If \(P\) is any numbered sequence of patterns, then \(Pv_i\)
is obtained from \(P\) by dividing rectangle \(i\) of \(P\)
vertically, giving the number \(i\) to the left half and the number
\(i+1\) to the right half, preserving the numbers of all rectangles
numbered less than \(i\) in \(P\), adding 1 to the numbers of all
rectangles numbered greater than \(i\) in \(P\).

Similar remarks apply to \(Ph_i\) except that rectangle \(i\) of
\(P\) is now divided horizontally with the lower half numbered \(i\)
and the upper half numbered \(i+1\).  Lastly \(P\sigma_i\) is
obtained from \(P\) by exchanging the numbers of rectangles numbered
\(i\) and \(i+1\) and making no other changes.

It is clear from the description of the multiplication in \(\Pi\) as
described in Section \ref{HatTwoVMonMult} that any sequence
of patterns (with perhaps the wrong numbering) can be obtained as a
word in the \(v_i\) and \(h_i\).  Now the numbering can be fixed up
by following this word with a word in the \(\sigma_i\).  Thus \(\Pi\)
is generated by \(\{v_i, h_i, \sigma_i\mid i\in \N\}\).

The previous paragraph says more.  It implies that every element of
\(\Pi\) can be written as a word in the \(v_i\) and \(h_i\) followed
by a word in the \(\sigma_i\).  This will be supported by the claims
that we make next.

\subsection{Relations in the positive monoid of
\protect\(\protect\widehat{2V}\protect\)}

We give four sets of relations that the elements \(v_i\), \(h_i\)
and \(\sigma_i\) satisfy.  The fact that they hold is important and
the reader should verify that they do so by drawing pictures.

The first set is 
\mymargin{ThompRels}
\begin{alignat}{2}
\label{ThompRels} 
x_jy_i &= y_ix_{j+1}, \qquad&&i<j, 
\end{alignat}
where the symbols \(x\) and
\(y\) come independently from the set of symbols \(\{h,v\}\).  We
refer to the relations in \tref{ThompRels} as the ``Thompson
relations'' because of their resemblance to the relations in
Thompson's group \(F\) and because of their power in reducing
infinite presentations to finite presentations.

The next set is 
\mymargin{PermRels}
\begin{alignat}{2}
\label{PermRels}
\sigma_i^2
&= 1, \qquad && i\ge0, \\ \sigma_i\sigma_j &= \sigma_j\sigma_i,
\qquad && |i-j|\ge2, \\ \sigma_i\sigma_{i+1}\sigma_i &=
\sigma_{i+1}\sigma_i \sigma_{i+1}, \qquad && i\ge0, 
\end{alignat} 
which are simply the relations of the infinite permutation group
generated by the transpositions \(\sigma_i\).

The third set is
\mymargin{ZSRels}
\begin{equation}
\label{ZSRels}
\sigma_j{x_i} =
\begin{cases}x_i \sigma_{j+1}, \quad&i<j, \\ x_{j+1}
\sigma_j\sigma_{j+1},
\quad& i=j, \\ x_j \sigma_{j+1}\sigma_j, \quad &i=j+1, \\ x_i \sigma_j,
\quad& i>j+1  \end{cases}
\end{equation}
where the symbol \(x\) comes from \(\{v,h\}\).  These relations give
the interaction between the \(v_i\) and \(h_i\) on the one hand  and
the \(\sigma_i\) on the other.  The relations \tref{ZSRels} tell
how to ``switch'' an \(x_i\) and \(\sigma_j\) that occur in the
``wrong'' order.  The relations \tref{ZSRels} give another argument for 
the claim that any element of \(\Pi\) can be written as a word in the
\(v_i\) and \(h_i\) followed by a word in the \(\sigma_i\).

The last set is 
\mymargin{CrossRels}
\begin{alignat}{2}
\label{CrossRels}
v_ih_{i+1}h_i &= h_iv_{i+1}v_i \sigma_{i+1}, \qquad&& i\ge0,
\end{alignat} which we refer to as the ``cross relations.''  A
picture is warranted to explain why.  On the left we show square 0
of \(v_0h_{1}h_0\) and on the right we show square 0 of
\(h_0v_{1}v_0\).  (Note the omission of \(\sigma_{1}\) from the
second expression.)

\[
\xy
(-9,-9); (-9,9)**@{-}; (9,9)**@{-}; (9,-9)**@{-}; (-9,-9)**@{-};
(-9,0); (9,0)**@{-};
(0,-9); (0,9)**@{-};
(-4.5,-4.5)*{0};
(4.5,-4.5)*{2};
(-4.5,4.5)*{1};
(4.5,4.5)*{3};
\endxy
\qquad
\qquad
\qquad
\xy
(-9,-9); (-9,9)**@{-}; (9,9)**@{-}; (9,-9)**@{-}; (-9,-9)**@{-};
(-9,0); (9,0)**@{-};
(0,-9); (0,9)**@{-};
(-4.5,-4.5)*{0};
(4.5,-4.5)*{1};
(-4.5,4.5)*{2};
(4.5,4.5)*{3};
\endxy
\]

\subsection{Remark}

There seems to be no group \(2F\).  The construction of the group
\(F\) uses the natural left-right order in the unit interval to
define a subgroup of \(V=1V\) of elements that preserve this order.
The author's attempts to select a preferred order for rectangles in
the unit square have resulted in failure.  The effort was to
construct a submonoid of \(\Pi\) that was cancellative, had common
right multiples and that respected some selected order of
rectangles.  The most successful attempt achieved cancellativity and
common right multiples, but was not a submonoid of \(\Pi\) since the
resulting multiplication was not associative.

\section{Generators for \protect\(2V\protect\)}

In this section, we find generators for \(2V\) and we will prove
that they generate.  In the next section, we will find relations for
\(2V\) but we will not prove that they suffice to give a
presentation.  However, in order to do what we want in this section
and the next, we need to be able to multiply elements.  We only have
to multiply in practice and not in theory.  Thus it will be
sufficient to give a method for multiplying elements that works for
all the calculations that we need, but we will not have to prove
that our method always works.  In fact it does, and the curious
reader can supply the reasons why.

Before we develop a method for mutliplying elements, we will work
out a nice method for representing elements that will make use of
our knowledge of the positive monoid \(\Pi\) of \(\widehat{2V}\).

\subsection{Elements as pairs of words}

Each element of \(2V\) is a homeomorphism from \(C\times C\) to
itself.  Each element of \(\widehat{2V}\) is a homeomorphism from
\(Y\), a countable union of copies of \(C\times C\), to itself.  We
can think of \(2V\) as a subgroup of \(\widehat{2V}\) by thinking of
\(2V\) as those elements of \(\widehat{2V}\) that take the copy of
\(C\times C\) in square 0 onto itself and that is the identity on
all other copies of \(C\times C\).

If we view elements of \(\widehat{2V}\) as pairs of numbered
sequences of patterns, then the elements of \(2V\) are those pairs
in which the patterns are both trivial after square 0, in which the
numbering in both patterns of the rectangles is consecutive after
square 0, and in which the count of rectangles in square 0 is the
same for both patterns.  Let \(\Pi_0\) be the set (it is a submonoid,
but in a very uninteresting way) of those elements of \(\Pi\)
corresponding to numbered sequences of patterns that are trivial
after square 0 and which number the squares consecutively after
square 0.  Then \(2V\) is given by pairs \((a,b)\) with \(a\) and
\(b\) in \(\Pi_0\) and so that the count of rectangles in square 0 is
the same for both \(a\) and \(b\).

Let \(w\) be a word in the elements of Section \ref{MonoidGens} in
the form discussed at the end of Section \ref{MonoidGens}.  That is,
\(w\) is a word in \(\{v_i, h_i, \sigma_i\mid i\in\N\}\) consisting
of a word \(p\) in \(\{v_i,h_i\mid i\in\N\}\) followed by a word
\(q\) in the \(\sigma_i\).  From the discussion in Section
\ref{MonoidGens}, we give sufficient conditions for such a word to
give and element of \(\Pi_0\).

We can build up the word \(p\) letter by letter from the trivial
word.  As we do so, we create approximations to the pattern
determined by \(p\).  The addition of each letter divides one
rectangle in the pattern built to that stage.  We assume that each
stage is a pattern in \(\Pi_0\).  Because of the nature of multiplying
on the right by a letter in \(\{v_i,h_i\mid i\in\N\}\), if the next
letter subdivides a rectangle in square 0, the result of adding that
letter will also give an element of \(\Pi_0\).  The trivial pattern is
in \(\Pi_0\), thus the word \(p\) gives an element of \(\Pi_0\) if
adding each letter to the prefix before it represents a division of
a rectangle in square 0.  It follows from this that if \[p=
x_{i_0}x_{i_1}\cdots x_{i_{k-1}}\] where each \(x\) independently
represents either the letter \(v\) or the letter \(h\) and \(i_j\le
j\) for \(0\le j < k\), then \(p\) represents an element of
\(\Pi_0\).

The letters in \(q\) are transpositions.  If \(p\) is as given
above, then the pattern that \(p\) represents has \(k+1\) rectangles
in square 0 numbered from \(0\) through \(k\).  As long as the
permutation given by \(q\) affects only rectangles 0 through \(k\),
the word \(w=pq\) will give an element in \(\Pi_0\).  Thus the
sufficient condition that we get is as follows.

\begin{lemma}\mylabel{InPiNought} Let \(w=pq\) be a word with \(p=
x_{i_0}x_{i_1}\cdots x_{i_{k-1}}\) a word in \(\{v_i,h_i\mid
i\in\N\}\) in that each \(x\) independently stands for \(v\) or
\(h\), and \(q=\sigma_{m_0}\sigma_{m_1}\cdots \sigma_{m_{n-1}}\) a
word in the \(\sigma_i\).  If \(i_j\le j\) for \(0\le j<k\) and
\(m_j < k\) for \(0\le j <n\), then \(w\) represents an element of
\(\Pi_0\).  \end{lemma}

We need a partial converse to Lemma \ref{InPiNought}.

\begin{lemma}\mylabel{InPiNoughtII} Let \(a\) be in \(\Pi_0\).  Then
there is a word \(w=pq\) satsifying the description in Lemma
\ref{InPiNought} that represents \(a\).  \end{lemma}

\begin{proof} We know that there is a word \(pq'\) with \(p=
x_{i_0}x_{i_1}\cdots x_{i_{k-1}}\) a word in \(\{v_i,h_i\mid
i\in\N\}\) and \(q'\) a word in the \(\sigma_i\) representing \(a\).
Part of the conclusion is forced.  If any \(i_j\) is greater than
\(j\), then some square other than square 0 is divided into smaller
rectangles.  This is something that cannot be undone, so we must
have \(i_j\le j\) for all \(j\).  Thus \(p\) represents an element
of \(\Pi_0\).  This means that the ordering of all the squares other
than square 0 is already correct.  The only way that this can still
be true of \(pq'\) is for \(q'\) to be the trivial permutation above
\(k\).  This means that \(q'\) is a permutation on \(\{0,1,\ldots,
k\}\) and can be rewritten as a word \(q\) in \(\{\sigma_0,
\sigma_1, \ldots, \sigma_{k-1}\}\).  The desired word is \(pq\).
\end{proof}

\subsection{Multiplying pairs of words}\mylabel{WordPairMult}

If \((a,b)\) is an element of \(\widehat{2V}\) with \(a\) and \(b\)
in \(\Pi\), then it represents the element \(ab^{-1}\) in
\(\widehat{2V}\).  Thus we get \((a,b)^{-1}=(b,a)\) and
\((a,b)(b,c)=(a,c)\).  This makes it easy to invert pairs \((a,b)\)
representing elements of \(2V\) with \(a\) and \(b\) in \(\Pi_0\)
since \((b,a)\) is still a pair of words with both entries in
\(\Pi_0\).

If \((a,b)\) and \((c,d)\) are pairs with all entries in \(\Pi_0\),
then we take advantange of the fact that \((ax,bx)\) represents the
same element as \((a,b)\) for any \(x\).  We only need to find \(x\)
and \(y\) so that \(bx=cy\) and so that all of \(ax\), \(bx\),
\(cy\) and \(dy\) are in \(\Pi_0\) to write \[(a,b)(c,d) =
(ax,bx)(cy,dy) = (ax,dy).\] In fact, this can always be done, but we
do not have to prove that.  As long as it always happens for the
calculations that we do, we are in good shape.

\subsection{The generators}\mylabel{TwoVeeGens}

The reader can verify that all of the following represent elements
of \(2V\):
\begin{alignat*}{2}
A_i &=(v_0^{i+1}v_1, \, v_0^{i+2}), &\qquad&i\ge0, \\
B_i &=(v_0^{i+1}h_1, \, v_0^{i+2}), &&i\ge0, \\
C_i &=(v_0^ih_0,\, v_0^{i+1}), &&i\ge0, \\
\pi_i &= (v_0^{i+2}\sigma_1, \, v_0^{i+2}), && i\ge0, \\
\opi_i &= (v_0^{i+1}\sigma_0, \, v_0^{i+1}), &&i\ge0.
\end{alignat*}  We let \[\Sigma = \{A_i, B_i, C_i, \pi_i, \opi_i\mid
i\in \N\}.\]

In the above list, the reader will note that \(C_0\) is the element
we called the baker's map in Section \ref{HomeoFromPatt}.  The
reader will also note that all second entries in the pairs are
powers of \(v_0\).  We let \(2V_0\) be the set of elements that can
be represented by a pair \((a,v_0^k)\) with \(a\) in \(\Pi_0\).
Note that \(k\) is determined from \(a\) since \(v_0^k\) cuts square
0 into \(k+1\) rectangles and \(a\) will only cut square 0 into
\(k+1\) rectangles if there are exactly \(k\) appearances of the
letters \(h\) or \(v\) in \(a\).

Now if \((a, v_0^k)\) and \((b, v_0^k)\) are elements of \(2V_0\)
(the fact that they use the same power if \(v_0\) is deliberate),
then \((a, v_0^k)(b, v_0^k)^{-1}=(a, v_0^k)(v_0^k,b)=(a,b)\) is an
element of \(2V\).  Conversely, if \((a,b)\) is in \(2V\), then
\(a\) and \(b\) divide square 0 into the same number (say \(k+1\)) of
rectangles and \((a,v_0^k)\) and \((b,v_0^k)\) are in \(2V\) with
\((a,b)=(a,v_0^k)(b,v_0^k)^{-1}\).  Thus if for every \(a\) in
\(\Pi_0\) we can show that \((a,v_0^k)\), with \(k+1\) the number of
rectangles in square 0 of the pattern for \(a\), is a word in the
elments of \(\Sigma\), then we will have shown that \(\Sigma\) is a
generating set for \(2V\).

In order to do this, we must understand how products of
elements in \(\Sigma\) behave.

\subsection{Building a pattern}

Let \(w=pq\) be a word as described in Lemma \ref{InPiNought}
representing an element of \(\Pi_0\).  We want to show that \((w,
v_0^k)\) can be obtained as a word in \(\Sigma\).  

We will modify some of the descriptions in Section \ref{TwoVeeGens}
to make them more convenient to work with.  To do this we take
advantage of the fact that for any \(c\), the pair \((ac,bc)\)
represents the same element \((a,b)\).

We give two sets of alternate formulations of some of the
generators.  
The second set of formulations is just re-indexing of the first.
The first set uses the relations \tref{ThompRels} and \tref{ZSRels}
to replace \(v_1v_0^{j-1}\) by \(v_0^{j-1}v_j\) and so forth.
\begin{alignat*}{2} A_i &=(v_0^{i+j}v_j,\, v_0^{i+j+1}),
&\qquad&i\ge0,\, j\ge1, \\ B_i &=(v_0^{i+j}h_j,\, v_0^{i+j+1}),
&&i\ge0,\,j\ge1,\\ \pi_i &= (v_0^{i+j+1}\sigma_j,\, v_0^{i+j+1}), &&
i\ge0,\,j\ge1, \\ \\ A_i &=(v_0^kv_{k-i},\, v_0^{k+1}), &&
i\ge0,\,k>i, \\ B_i &=(v_0^kh_{k-i},\, v_0^{k+1}), &&i\ge0,\,k>i,\\
\pi_i &= (v_0^{k+1}\sigma_{k-i},\, v_0^{k+1}), &&i\ge0,\, k>i.
\end{alignat*}

We are now ready to build \((w, v_0^k)\) with \(w=pq\) as in Lemma
\ref{InPiNought}.  We start with \((p, v_0^k)\).  Lemma
\ref{InPiNought} implies that if \(p\) is the empty word, then so is
\(q\), so we assume that \(p\) has at least one letter.

We have \(p= x_{i_0}x_{i_1}\cdots x_{i_{k-1}}\) a word in
\(\{v_i,h_i\mid i\in\N\}\) with each \(i_j\le j\).  This makes
\(x_{i_0}\) equal to \(v_0\) or \(h_0\).  We will build \((p,
v_0^k)\) as a word in \(\Sigma\), and the first letter of this word
is the trivial element represented as \((v_0, v_0)\) if
\(x_{i_0}=v_0\), and the first letter is \(C_0=(h_0,v_0)\) if
\(x_{i_0}=h_0\).  

Now for \(j\) with \(0<j<k-1\), we let \(p_j\) be the prefix
\(x_{i_0}x_{i_1}\cdots x_{i_{j-1}}\) of \(p\) of length \(j\) and
assume that \((p_j, v_0^j)\) has been represented as a word in
\(\Sigma\).  There are four cases to consider: \[x_{i_{j}} =
\begin{cases} v_0, \\ v_m, &0<m\le j, \\ h_0, \\ h_m, & 0<m\le
j. \end{cases}\] In the first case, we multiply by the identity
since \((p_j, v_0^j)\) represents the same element as \((p_jv_0,
v_0^{j+1}) = (p_{j+1}, v_0^{j+1})\).  In the second case, we write
\[(p_{j+1}, v_0^{j+1}) = (p_jv_m, v_0^{j+1}) = (p_jv_m,
v_0^jv_m)(v_0^jv_m, v_0^{j+1}) = (p_j, v_0^j) A_{j-m}.\] Similarly,
the third case leads to \((p_{j+1}, v_0^{j+1}) = (p_j, v_0^j)C_j\)
and the last case leads to \((p_{j+1}, v_0^{j+1}) = (p_j,
v_0^j)B_{j-m}\).  Inductively, we get to \((p,v_0^k)\) as a word in
\(\Sigma\).

Now we know that \(q\) is a word in \(\{\sigma_0, \sigma_1, \ldots,
\sigma_{k-1}\}\).  We have available in \(\Sigma\) the elements
\(\opi_{k-1} = (v_0^k\sigma_0, v_0^k)\) and
\(\pi_i=(v_0^k\sigma_{k-1-i}, v_0^k)\) for \(0\le i<k-1\).  Setting
\(m=k-1-i\), gives \(0<m\le k-1\) and \(\pi_{k-1-m} = (v_0^k
\sigma_m, v_0^k)\).  As three of the four cases just above were
handled, it is now elementary that \((w,v_0^k) = (pq, v_0^k)\) can
be written as a word in \(\Sigma\).

From Lemma \ref{InPiNoughtII}, we know that any element of \(\Pi_0\)
is represented by a word as described in Lemma \ref{InPiNought}.
Thus we have shown the following.

\begin{prop}\mylabel{TwoVeeGensWork} The set \(\Sigma\) is a
generating set for the group \(2V\).  \end{prop}

\section{Relations for \protect\(2V\protect\)}

There is no subtlety in this section.  We have the following.

\begin{prop}\mylabel{TwoVeeRels}  The following 17 infinite
families of relations hold in \(2V\).  In the following the letters
\(X\) and \(Y\) represent symbols from \(\{A,B\}\).
\begin{alignat*}{2} 
X_qY_m &= Y_mX_{q+1},&\quad&m<q, \\
\pi_qX_m &= X_m\pi_{q+1}, &&m<q, \\
\pi_qX_q &= X_{q+1}\pi_q\pi_{q+1}, &&q\ge0, \\
%\pi_qX_{q+1} &= X_q\pi_{q+1}\pi_q, &&q\ge0, \\
\pi_qX_m &= X_m\pi_q, && m>q+1, \\
\opi_qX_m &= X_m\opi_{q+1}, &&m<q, \\
\opi_mA_m &= \pi_m\opi_{m+1}, &&m\ge0, \\
\opi_mB_m &= C_{m+1} \pi_m \opi_{m+1}, &&m\ge0, \\
%\opi_m C_{m+1} &= B_m \opi_{m+1}\pi_m, && m\ge0, \\
C_qX_m &= X_mC_{q+1}, &&m<q, \\
C_mA_m &= B_mC_{m+2}\pi_{m+1}, &&m\ge0, \\
\pi_qC_m &= C_m\pi_q, && m>q+1, \\
A_mB_{m+1}B_m &= B_mA_{m+1}A_m \pi_{m+1}, &&m\ge0, \\
\pi_q\pi_m &= \pi_m\pi_q, 
  &&|m-q|\ge2, \\
\pi_m\pi_{m+1}\pi_m &= \pi_{m+1}\pi_m\pi_{m+1},
  &&m\ge0, \\
\opi_q\pi_m &= \pi_m\opi_q, 
  &&q\ge m+2, \\
\pi_m\opi_{m+1}\pi_m &= \opi_{m+1}\pi_m\opi_{m+1},
  &&m\ge0, \\
\pi_m^2 &= 1, &&m\ge0, \\
\opi_m^2 &= 1, &&m\ge0.
\end{alignat*}
\end{prop}

It will be shown in \cite{brin:hd4} that the above list of relations
suffices to present \(2V\).

\begin{proof} The above is nothing but calculation.  We illustrate
some and leave the rest for the reader.  

For the first line with \(X=Y=A\), we assume \(m<q<k\).  
\[\begin{split} A_qA_m &=
(v_0^kv_{k-q},\,v_0^{k+1}) (v_0^{k+1}v_{k+1-m},\, v_0^{k+2}) \\ &=
(v_0^kv_{k-q}v_{k+1-m},\, v_0^{k+2}) \\ &= (v_0^kv_{k-m}v_{k-q}, \,
v_0^{k+2}) \\ &= (v_0^kv_{k-m},\, v_0^{k+1})
(v_0^{k+1}v_{k+1-(q+1)},\, v_0^{k+2}) \\ &=
A_mA_{q+1}. \end{split}\]

Similarly, we get all of \[\begin{split} A_qB_m &= B_mA_{q+1}, \\
B_qB_m &= B_mB_{q+1}, \\ B_qA_m &= A_mB_{q+1}. \end{split}\]

We now turn to the \(\pi_i\).  We note that \begin{alignat*}{3}
m&<q &\quad&\rightarrow \quad& (k+1-m) &> (k-q)+1, \\ 
m&=q &&\rightarrow & (k+1-m) &= (k-q)+1, \\ 
m&=q+1 &&\rightarrow & (k+1-m) &= (k-q), \\ 
m&>q+1 &&\rightarrow & (k+1-m) &< (k-q).
\end{alignat*}

Now we get \[\begin{split} \pi_qA_m &=(v_0^{k+1}\sigma_{k-q},\,
v_0^{k+1})(v_0^{k+1} v_{k+1-m},\, v_0^{k+2}) \\ &= (v_0^{k+1}\sigma_{k-q}
v_{k+1-m}, \, v_0^{k+2}) \\ &= \begin{cases}
(v_0^{k+1}v_{k+1-m}\sigma_{k-q},\, v_0^{k+2}) = A_m\pi_{q+1}, &\quad m<q,
\\ \begin{split}
(v_0^{k+1}v_{k-m}\sigma_{k+1-m}\sigma_{k-m},\, v_0^{k+2}) &=
A_{m+1}\pi_m\pi_{m+1} \\  &= A_{q+1}\pi_q\pi_{q+1},\end{split}
 &\quad m=q, \\ \begin{split}
(v_0v_{k+2-m}\sigma_{k-q}\sigma_{k+1-q},\, v_0^{k+2}) &=
A_{m-1}\pi_{q+1}\pi_{q}\\ &= A_q\pi_{q+1}\pi_{q},\end{split}
 &\quad m=q+1, \\
(v_0^{k+1}v_{k+1-m}\sigma_{k+1-q},\, v_0^{k+2}) = A_m\pi_q, & \quad
m>q+1. \end{cases}
\end{split}\]
Note that the case \(m=q+1\) (not mentioned in the statement of the
proposition) actually follows from the case \(m=q\) since each
\(\pi_i\) is its own inverse.

We get similar calculations for \(\pi_qB_m\).

Now we look at \[\begin{split} \opi_mB_m &= (v_0^{m+1}\sigma_0,\,
v_0^{m+1})(v_0^{m+1} h_1,\, v_0^{m+2}) \\ &= (v_0^{m+1}
\sigma_0h_1,\, v_0^{m+2}) \\ &= (v_0^{m+1} h_0 \sigma_1 \sigma_0,\,
v_0^{m+2}) \\ &= (v_0^{m+1}h_0,\, v_0^{m+2})(v_0^{m+2} \sigma_1,\,
v_0^{m+2})(v_0^{m+2} \sigma_0,\, v_0^{m+2}) \\ &= C_{m+1} \pi_m
\opi_{m+1}. \end{split}\]

One last example is \[\begin{split} C_mA_m &= (v_0^mh_0v_0,\,
v_0^{m+2}) (v_0^{m+2}v_2,\, v_0^{m+3}) \\ &= (v_0^mh_0v_0v_2,\,
v_0^{m+3}) \\ &= (v_0^mh_0v_1v_0,\, v_0^{m+3}) \\ &=
(v_0^mv_0h_1h_0\sigma_1, v_0^{m+3}) \\ &= (v_0^{m+1}h_1,\,
v_0^{m+2})(v_0^{m+2}h_0,\, v_0^{m+3})(v_0^{m+3}\sigma_1,\,
v_0^{m+3}) \\ &= B_m C_{m+2} \pi_{m+1}. \end{split}\] 

All of the above are supported by drawing pictures.  \end{proof}

An immediate consequence of Proposition \ref{TwoVeeRels} is that
\(2V\) is finitely generated.

\begin{prop}\mylabel{TwoVeeFiniteGens} The group \(2V\) is generated
by \(\{A_i, B_i, \pi_i, \opi_i\mid i=0,1\}\).  \end{prop}

\begin{proof} The relations in Proposition \ref{TwoVeeRels} include
relations of the type \(Z_qA_m=A_mZ_{q+1}\) whenever \(m<q\) and
\(Z\) is from \(\{A,B,C,\pi,\opi\}\).  From this
\(Z_{q+1}=A_0^{-q}Z_1A_0^q\) follows and we see that \(2V\) is
generated by \(\{A_i, B_i, C_i, \pi_i, \opi_i\mid i=0,1\}\).

From \(C_mA_m = B_mC_{m+2}\pi_{m+1}\) and \(C_{m+1}B_m=B_mC_{m+2}\)
we get \(C_mA_m = C_{m+1}B_m\pi_{m+1}\) which gives \(C_m =
C_{m+1}(B_m\pi_{m+1}A_m^{-1})\).  Now we use \(\opi_mB_m = C_{m+1}
\pi_m \opi_{m+1}\) to get \(C_{m+1} = \opi_mB_m\opi_{m+1}\pi_m\)
using the fact that the \(\pi_i\) and \(\opi_i\) are involutions.
Combining the two gives \(C_m =
(\opi_mB_m\opi_{m+1}\pi_m)(B_m\pi_{m+1}A_m^{-1})\).  \end{proof}

\section{Simplicity of \protect\(2V\protect\)}

From Proposition \ref{SimplComm}, we must show the following.

\begin{prop}\mylabel{TwoVeeAbel} The group \(2V\) equals its
commutator subgroup.  \end{prop}

\begin{proof} This is a very straightforward calculation showing
that the eight generators of \(2V\) are all products of commutators.
To eliminate a few words, we use \(V\simeq W\) to mean that \(V=W\)
modulo the commutator subgroup.

From the proof of Proposition \ref{TwoVeeFiniteGens}, we already
know that \(Z_q\simeq Z_1\) for all \(q>1\) and \(Z\) in
\(\{A,B,C,\pi,\opi\}\).

From \(\pi_0A_0=A_1\pi_1\pi_0\), we get \(A_0\simeq A_1\pi_1\).
Similarly, \(B_0\simeq B_1\pi_1\).

Now \(\pi_1A_1=A_2\pi_1\pi_2\) gives \(\pi_1A_1\simeq
A_1\pi_1\pi_1=A_1\) and \(\pi_1\simeq 1\).  Thus \(A_0\simeq A_1\)
and \(B_0\simeq B_1\).  Also \(\pi_0\pi_1\pi_0=\pi_1\pi_0\pi_1\)
gives \(\pi_0\simeq\pi_1\), so \(\pi_0\simeq1\).

From \(\opi_1A_1=\pi_1\opi_2\simeq \opi_1\) we get \(A_0\simeq
A_1\simeq 1\).

The pair of relations \(\opi_0A_0=\pi_0\opi_1\) and
\(\pi_0\opi_1\pi_0=\opi_1\pi_0\opi_1\) give \(\opi_0\simeq
\opi_1\simeq 1\).

Lastly \(A_0B_1B_0=B_0A_1A_0\pi_1\) gives \(B_1\simeq 1\) so
\(B_0\simeq B_1\simeq 1\).  This shows that all 8 generators from
Proposition \ref{TwoVeeFiniteGens} are in the commutator subgroup.
\end{proof}

From Propositions \ref{SimplComm} and \ref{TwoVeeAbel} we get the
following.

\begin{thm}\mylabel{TwoVeeSimpl} The group \(2V\) is simple.
\end{thm}

\section{On the baker's map}

The baker's map given by 

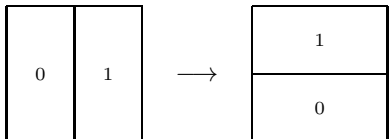
\[ \xy (-9,-9); (-9,9)**@{-}; (9,9)**@{-}; (9,-9)**@{-};
(-9,-9)**@{-}; (0,-9); (0,9)**@{-}; (-4.5,0)*{\scriptstyle0};
(4.5,0)*{\scriptstyle1}; \endxy \quad \longrightarrow \quad \xy
(-9,-9); (-9,9)**@{-}; (9,9)**@{-}; (9,-9)**@{-}; (-9,-9)**@{-};
(-9,0); (9,0)**@{-}; (0,-4.5)*{\scriptstyle0};
(0,4.5)*{\scriptstyle1}; \endxy \] distinguishes \(2V\) sharply from
\(V\).  The following lemma and corollary are well known.  We
learned of them from Dennis Pixton.

Let elements of \(C\times C\) be represented by based doubly
infinite strings in \(\{0,1\}\) in the following manner.  The first
coordinate in \(C\times C\) will be represented by an infinite
sequence of elements of \(\{0,1\}\) written from left to right.  The
second coordinate in \(C\times C\) will be represented by an
infinite sequence of elements of \(\{0,1\}\) written from right to
left.  The two sequences are then written on the same line separated
by a ``binary point.''  Obviously, this can be viewed as a function
from the integers \(\Z\) to \(\{0,1\}\).  If this function is
written so that \(x_i\) is the image of \(i\), then the binary point
can be viewed as coming between \(x_{-1}\) and \(x_0\) so that the
sequence \((x_i)_{i\ge0}\) gives the first coordinate in \(C\times
C\) and \((x_i)_{i<0}\) gives the second coordinate.

\begin{lemma}\mylabel{BakerShift} The baker's map corresponds to
shifting a based doubly infinite sequence from \(\{0,1\}\) one
position.  Specifically, if \(b\) is the baker's map and
\(x:\Z\rightarrow \{0,1\}\) is a sequene representing a point in
\(C\times C\), then \((b(x))_i=x_{i+1}\).  \end{lemma}

\begin{proof} There are two cases to consider.  If the first
coordinate starts with 0, then the point is from the left half of
the unit square in the figure above, and if the first coordinate
starts with 1, then the point is from the right half of the unit
square.  The remaining details are easy.  \end{proof}

\begin{cor}\mylabel{UnbddOrbits} There is no bound on the size of
the finite orbits of the baker's map.  \end{cor}

\begin{proof} A periodic function \(x:\Z\rightarrow \{0,1\}\) with
period \(p\) lies in a finite orbit of the baker's map of size
\(p\).  \end{proof}

\section{Rubin's theorem}\mylabel{RubinThmSec}

The previous section shows that the baker's map has complex
dynamics.  In the next section we will show that no element of \(V\)
or a related group has such dynamics.  We will make use of these
observations by applying results of Mati Rubin.  Generally, these
say that under a set of hypotheses, the groups involved in two group
actions are isomorphic if and only the actions are topologically
conjugate.  Since our observations will show that the actions of
\(V\) on \(C\) and \(2V\) on \(C\times C\) cannot be topologically
conjugate, we will have that \(V\) and \(2V\) are not isomorphic.

\subsection{The theorem} We now give the definitions needed to state
Rubin's result.  If \(X\) is a topological space, if \(H(X)\) is its
group of self homeomorphisms, and if \(F\sub H(X)\) is a subgroup of
\(H(X)\), then we say that \(F\) is {\itshape locally dense} if for
every \(x\in X\) and every open \(U\sub X\) with \(x\in U\), the
closure of \[\{f(x)\mid f\in F,\,\,\,
f|_{(X-U)}=\mathbf{1}_{(X-U)}\}\] contains some open set.  In other
words, for each open \(U\), the subgroup of elements fixed off \(U\)
has every orbit in \(U\) dense in some open set in \(U\).

The following is essentially Theorem 3.1 of \cite{MR99d:20003} where
it is described as a combination of parts (a), (b) and (c) of
Theorem 3.5 of \cite{MR90g:54044}.  The hypothesis that there be no
isolated points was inadvertently omitted from \cite{MR99d:20003}
where it is needed.  The terminology \emph{locally dense} is not
used in either \cite{MR99d:20003} or \cite{MR90g:54044}.  However,
in the absence of isolated points, it implies the notion of
\emph{locally moving} that is used in \cite{MR99d:20003}.  The
absense of isolated points seems to correspond to the assumption of
``no atoms'' in the Boolean algebras of \cite{MR90g:54044}.

\begin{thm}[Rubin]\mylabel{RubinThm} Let \(X\) and \(Y\) be locally
compact, Hausdorff topological spaces without isolated points, let
\(H(X)\) and \(H(Y)\) be the self homeomorphism groups of \(X\) and
\(Y\), respectively, and let \(G\sub H(X)\) and \(H\sub H(Y)\) be
subgroups.  If \(G\) and \(H\) are isomorphic and are both locally
dense, then for each isomorphism \(\phi:G\rightarrow H\) there is a
unique homeomorphism \(\tau:X\rightarrow Y\) so that for each \(g\in
G\), we have \(\phi(g) = \tau g \tau^{-1}\).  \end{thm}

\subsection{The scope of Rubin's theorem}\mylabel{RubinScopeSec}The
Cantor set \(C\) is a locally compact, Hausdorff topological space
with no isolated points as is the homeomorphic \(C\times C\).  That
the action of \(2V\) on \(C\times C\) is locally dense follows
immediately from the construction of the elements.

We will want to apply Rubin's theorem to other simple groups.
Groups \(V_{n,r}\) are introduced in \cite{higman:fpsimpl} (where
they are called \(G_{n,r}\)) that are generalizations of what we
call \(1V\).  Patterns are created in \(r\) disjoint copies of the
unit interval and these are used to create self homeomorphisms of
\(r\) disjoint copies of the Cantor set.  The first pattern consists
of the \(r\) separate, original unit intervals and new pattens are
created by subdividing intervals into \(n\) equal subintervals.
Thus \(1V=V_{2,1}\).  

It is shown in \cite{higman:fpsimpl} that the commutator subgroups
\(V^+_{n,r}\) of all the \(V_{n,r}\) are infinite, simple, and
finitely presented.  Further, in \cite{brown:finiteprop}, subgroups
\(T^s_{n,r}\) of the \(V_{n,r}\) (again, the \(V_{n,r}\) are called
\(G_{n,r}\) in \cite{brown:finiteprop}) are given and are shown to
be infinite, simple, and finitely presented.  The groups
\(T^s_{n,r}\) are the second commutator subgroups of larger groups
\(T_{n,r}\) that are also studied in \cite{brown:finiteprop}.  It is
known that Proposition \ref{SmallSupsGen} applies to the \(V_{n,r}\)
and \(T_{n,r}\) as well as a strengthening of Lemma \ref{SuffTrans}
that says that if \(K\) and \(U\) of that lemma are contained in an
open \(V\), then the element \(h\) of that lemma can be chosen to
have its support in \(V\).  It is now easy to argue that the given
actions of the \(V_{n,r}\) and \(T_{n,r}\) are locally dense as well
as the restriction of that action to any member of the derived
series.  Thus Rubin's theorem applies to all of the \(V^+_{n,r}\)
and \(T^s_{n,r}\).

\section{The dynamics of elements of
\protect\(V\protect\)}\mylabel{VDymSec}

The purpose of this section is to prove the following.

\begin{prop}\mylabel{BddOrbits} Let \(f\) be an element of \(V\).
Then there is an \(n(f)\in \N\) so that any finite orbit of \(f\)
has no more than \(n(f)\) elements.  \end{prop}

Trivial modifications of the arguments give the same results for the
\(V_{n,r}\).  Since all of the groups mentioned after the statement
of Rubin's theorem are subgroups of the \(V_{n,r}\), we get similar
results for them.  From Proposition \ref{BddOrbits}, Corollary
\ref{UnbddOrbits} and Theorem \ref{RubinThm}, we get the following.

\begin{thm}\mylabel{NoIso} The group \(2V\) is isomorphic to none of
the infinite, simple, finitely presented groups \(V^+_{n,r}\) and
\(T^s_{n,r}\).  \end{thm}

\subsection{The scope of Theorem
\protect\ref{NoIso}}\mylabel{NoIsoScopeSec}The papers \cite{scott}
and \cite{MR2001g:20028} give constructions of other infinite,
simple, finitely presented groups.  We do not know if Theorem
\ref{NoIso} can be extended to cover these groups.  There are also
groups with these properties constructed in \cite{MR2002i:20042}.
The theorem does cover the groups constructed in
\cite{MR2002i:20042} for trivial reasons: the groups in
\cite{MR2002i:20042} are all torsion free.

\subsection{Strategy}

Elements of \(V\) can be determined by pairs of trees plus a
permutation.  The proof proceeds by  modifying a tree-pair
of an arbitrary element \(f\) of \(V\) until it is possible to read
all of the dynamics of \(f\) from the tree-pair.

We note that the modified tree-pair will not be the pair that is
usually thought of as giving a ``normal form'' for \(f\).  The
normal form is the smallest in size (see \cite{brin:bv}), while the
pair that reveals the dynamics is usually not the smallest.

\subsection{Elements as pairs of trees}

The statement \(V=1V\) should enable the reader to describe the
elements of \(V\): they are given by pairs of numbered patterns of
the unit interval \(I=[0,1]\).  The discussion is identical to that
in Section \ref{DefTwoVeeSec}.  We will alter the description to
pairs of trees with a permutation for two reasons.  It is
traditional (see \cite{CFP} and \cite{brown:finiteprop}) and we will
find it useful.

Let \(T\) be the set of finite words (including the empty word) on
\(\{0,1\}\).  It is a monoid under concatenation (and in fact the
free monoid on two generators) with the empty word \(\phi\) as the
identity.  We also think of it as the infinite binary tree (and we
refer to the elements of \(T\) as \emph{nodes} when we do) since we
can think of \(v0\) and \(v1\) as the (respectively, left and right)
child nodes of the node \(v\in T\).  The empty word \(\phi\) is the
root node of \(T\).

Each node in \(T\) corresponds to an interval in a pattern on
\(I\).  Recursively, \(\phi\) corresponds to \(I\) itself and if
\(v\) corresponds to \([a,b]\), then \(v0\) corresponds to \([a,c]\)
and \(v1\) corresponds to \([c,b]\) where \(c=(a+b)/2\).

For us, a finite tree will be a finite subset \(D\) of \(T\) so that
(1) every prefix of a node in \(D\) is also in \(D\), and (2) \(v0\)
is in \(D\) if and only if \(v1\) is in \(D\).  The \emph{leaves} of
such a \(D\) will be the nodes in \(D\) whose children are not in
\(D\).  Nodes of a tree that are not leaves are called
\emph{interior nodes} of a tree.  Below is a picture of a finite
tree with five leaves and four interior nodes.

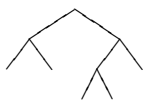
\[
\xy
(-4,0); (-1,4)**@{-}; (5,8)**@{-}; (11,4)**@{-}; (8,0)**@{-};
(6,-4)**@{-}; 
(-1,4); (2,0)**@{-};
(11,4); (14,0)**@{-};
(8,0); (10,-4)**@{-};
\endxy
\]

Note that the root of every tree is the empty word \(\phi\).

The leaves of a finite tree \(D\) give a pattern in \(I\) by taking
the intervals in \(I\) corresponding to the leaves of \(D\).  Two
trees \(D\) and \(R\) (for Domain and Range) with the same number
\(n\) of leaves define two patterns in \(I\) with the same number of
intervals.  If we are now given a one-to-one correspondence between
the intervals obtained from \(D\) to the intervals obtained from
\(R\), then we can build a homeomorphism from the Cantor set \(C\)
to itself in a manner analogous to that in Section
\ref{HomeoFromPatt}.

We now think of elements of \(V\) as triples \((D, \sigma, R)\)
where \(D\) and \(R\) are finite trees with the same number \(n\) of
leaves and where \(\sigma\) is a bijection from the leaves of \(D\)
to the leaves of \(R\).

Note that \(\sigma\) can be replaced by numberings of the leaves of
\(D\) and \(R\) (as is done in \cite{CFP}), but the triple notation
(as used in \cite{brown:finiteprop}) will be more convenient for us.

\subsection{Tree operations and carets}

Insisting that every finite tree be a subset of a single infinite
tree \(T\) will have its advantanges.  Given a triple \((D, \sigma,
R)\), we will have reason to refer to \(D\cap R\), to \(D-R\), and
to \(R-D\) which are now nicely defined.  We improve on the niceness
by introducing carets.  A caret is any triple \((v, v0, v1)\) in
\(T\).  Every finite tree is a finite union of carets if we sloppily
declare that the trivial tree \(\{\phi\}\) is the union of zero
carets.

We will insist on this view when we write down \(D-R\) and \(R-D\)
and will say that \(D-R\) is the set of carets in \(D\) that are not
in \(R\).  It is seen that \(D-R\) breaks up into a union of
pairwise disjoint ``trees'' (called the \emph{components} of
\(D-R\)) whose roots are not the empty word \(\phi\), but are leaves
of \(R\).  Similar remarks apply to \(R-D\).

In the picture below are two trees with 5 leaves, one \(D\) with
solid lines and the other \(R\) with dashed lines.  They are drawn
slightly offset so they can both be seen.  In the picture, \(D-R\)
has one component with two carets, and \(R-D\) has two components
with one caret each.  It is important to keep in mind that \(D-R\)
and \(R-D\) are differences of sets of carets.

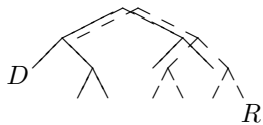
\[
\xy
(-2,-6); (-4,-2)**@{-}; (-8,2)**@{-}; (0,6)**@{-};
(8,2)**@{-}; (4,-2)**@{-};
(-4,-2); (-6,-6)**@{-}; (-8,2); (-12,-2)**@{-};
(8,2); (12,-2)**@{-};
(-6,2); (2,6)**@{--}; (10,2)**@{--}; (6,-2)**@{--}; (8,-6)**@{--};
(6,-2); (4,-6)**@{--}; (10,2); (14,-2)**@{--}; (12,-6)**@{--};
(14,-2); (16,-6)**@{--};
(-14,-3)*{D}; (17,-8)*{R};
\endxy
\]

An important triviality is that the number of leaves of a finite
tree is one more than the number of carets in the tree.  Since
number of carets is clearly the number of interior nodes, the number
of leaves is one more than the number of interior nodes.

\subsection{Altering tree pairs}

If \((D,\sigma,R)\) represents an element of \(V\), then we can
modify the triple to create another representing the same element.
Let \(u\) be one of the leaves of \(D\) and let \(U\) be a binary
tree with \(k\) leaves.  We can create a new triple \((D', \sigma',
R')\) from \((D, \sigma, R)\) and the pair \((u , U)\) as follows.

We attach a copy of \(U\) to \(D\) to form \(D'\) by identifying the
root of \(U\) to \(u \).  Note that this is simply forming the tree
\(D'=D\cup uU\) where \(uU\) is just the product (concatenation) of
\(u\) with the nodes of \(U\).  We get \(R'\) as \(R\cup
\sigma(u)U\).  Note that the leaves of \(D'\) that are not leaves of
\(D\) are of the form \(uv\) as \(v\) runs over the leaves of \(U\).
The leaves of \(R'\) that are not leaves of \(R\) are \(\sigma(u)v\)
as \(v\) runs over the leaves of \(U\).

We define \(\sigma'\) so that it agrees with \(\sigma\) on leaves
that \(D'\) shares with \(D\) and so that it takes each leaf of the
form \(uv\) with \(v\) a leaf of \(U\) to \(\sigma(u)v\).  We call
the triple \((D', \sigma', R')\) an \emph{augmentation} of
\((D,\sigma,R)\) at leaf \(u \) by \(U\).  The next figure gives an
element where the bijection is illustrated by the labels on the
leaves and the result of the modification is by the pair \((b,
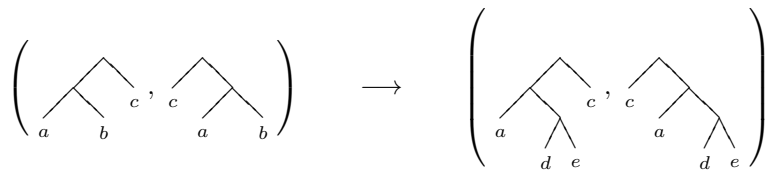
\xy(0,-1);(2,1)**@{-};(4,-1)**@{-}\endxy)\).

\[
\left(
\,\,
\xy
(-4,-4);(0,0)**@{-};(4,4)**@{-};(8,0)**@{-};
(0,0);(4,-4)**@{-};
(-4,-6)*{\scriptstyle a};
(4,-6)*{\scriptstyle b};
(8,-2)*{\scriptstyle c};
\endxy
\,\,,\,\,
\xy
(4,-4);(0,0)**@{-};(-4,4)**@{-};(-8,0)**@{-};
(0,0);(-4,-4)**@{-};
(4,-6)*{\scriptstyle b};
(-4,-6)*{\scriptstyle a};
(-8,-2)*{\scriptstyle c};
\endxy
\,\,
\right)
\qquad
\longrightarrow
\qquad
\left(
\,\,
\xy
(-4,-4);(0,0)**@{-};(4,4)**@{-};(8,0)**@{-};
(0,0);(4,-4)**@{-};
(2,-8);(4,-4)**@{-};(6,-8)**@{-};
(-4,-6)*{\scriptstyle a};
(2,-10)*{\scriptstyle d};
(6,-10)*{\scriptstyle e};
(8,-2)*{\scriptstyle c};
\endxy
\,\,,\,\,
\xy
(4,-4);(0,0)**@{-};(-4,4)**@{-};(-8,0)**@{-};
(0,0);(-4,-4)**@{-};
(2,-8);(4,-4)**@{-};(6,-8)**@{-};
(2,-10)*{\scriptstyle d};
(6,-10)*{\scriptstyle e};
(-4,-6)*{\scriptstyle a};
(-8,-2)*{\scriptstyle c};
\endxy
\,\,
\right)
\]

It is elementary that if \((D', \sigma', R')\) is obtained from
\((D, \sigma, R)\) by an augmentation, then \((D',\sigma',R')\)
represents the same element of \(V\) as \((D, \sigma,R)\).

\subsection{Iterated augmentations}

Augmentations can be iterated under certain conditions.  Assume we
have a sequence of leaves \(u_1, \ldots, u_n\) of \(D\) so that all
of \(u_1, \ldots, u_n\) and \(\sigma(u_n)\) are different and so
that \(u_{i+1}=\sigma(u_i)\) for \(1\le i <n\).  Here we take
advantage of the fact that all trees are rooted subtrees of the
complete binary tree \(T\).  Now given a tree \(U\), we can add, for
each \(1\le i\le n\), a copy of \(U\) to \(D\) at \(u_i\) and a copy
of \(U\) to \(R\) at \(\sigma(u_i)\).  There is a resulting one-to-one
correspondence that makes the resulting triple a representative of
the same element as represented by \((D, \sigma, R)\).  We call the
alterations to \(R\) and \(D\) just described as an {\itshape
iterated augmentation by \(U\) along} \(u_1, \ldots, u_n\) and we
call the sequence of leaves \(u_1, \ldots u_n\) that it is based on
an {\itshape iterated augmentation chain}.

\subsection{The argument}

We now look at particular properties of the triple \((D, \sigma,
R)\) that represents the element \(f\) of \(V\).

The number of carets of \(D\) equals the number of carets of \(R\)
since the number of leaves of \(D\) equals the number of leaves of
\(R\).  Let \(n\) be the number of carets of \(D\) and of \(R\), and
let \(m\) be the number of carets of \(D\cap R\).  Recall that
\(D-R\) refers to carets of \(D\) that are not in \(R\).  Then each
of \(D-R\) and \(R-D\) has \(n-m\) carets.

We call \(n-m\) the {\itshape imbalance} of the representative
\((D,\sigma,R)\) of \(f\).  We assume from this point on that
\((D,\sigma, R)\) has the least imbalance of all the representatives
of \(f\).

We now look at the number of components of \(D-R\) and \(R-D\).
Among the repsentatives with minimial imbalance, we take a
representative that has the smallest number of components of
\(D-R\).  Among all such representatives, we take a representative
that has the smallest number of components of \(R-D\).  We now prove
facts about any representative triple \((D, \sigma, R)\) chosen in
this manner.

\begin{lemma}\mylabel{FirstTechLemma} It is impossible to have an
iterated augmentation chain \(u_1, \ldots, u_n\) so that \(u_1\) is
an interior node of \(R\) and \(\sigma(u_n)\) is an interior node of
\(D\).  \end{lemma}

\begin{proof} If false, then with \(U\) the component of \(D-R\)
with root at \(\sigma(u_n)\), we can perform an iterated
augmentation by \(U\) along \(u_1, \ldots , u_n\).  This removes the
copy of \(U\) with root at \(\sigma(u_n)\) from \(D-R\).  There are
copies of \(U\) added to both \(D\) and \(R\) at each \(u_i\) with
\(2\le i\le n\) so these do not contribute to \(D-R\).  Now if \(n\)
is the number of carets in \(U\), then the number of carets from the
copy of \(U\) added to \(D\) at \(u_1\) that are added to \(D-R\) is
strictly fewer than \(n\) since \(u_1\) is an interior node of
\(R\).  This would lower the imbalance which is not possible by
choice.  \end{proof}

\begin{lemma} It is impossible to have an iterated augmentation
chain \(u_1,\ldots, u_n\) so that \(u_1\) is not a node of \(R\),
so that \(\sigma(u_n)\) is an interior node of \(D\) and so that
the component of \(D-R\) containing \(u_1\) is not the component of
\(D-R\) whose root is at \(\sigma(u_n)\).  \end{lemma}

\begin{proof} Let \(U\) be the component of \(D-R\) whose root is at
\(\sigma(u_n)\), and let \(V\) be the component of \(D-R\) that
contains the leaf \(u_1\) of \(D\).  An iterated augmentation by
\(U\) along \(u_1, \ldots, u_n\) would remove \(U\) as a component
of \(D-R\), would introduce no new components of \(D-R\) and would
add a copy of \(U\) to the component \(V\) of \(D-R\).  This
iterated augmentation would leave the imbalance unchanged, and would
reduce the number of components of \(D-R\).  This is not possible by
choice.  \end{proof}

\begin{lemma} It is impossible to have an iterated augmentation
chain \(u_1,\ldots, u_n\) so that \(\sigma(u_n)\) is not a node of
\(D\), so that \(u_1\) is an interior node of \(R\) and so that
the component of \(R-D\) containing \(\sigma(u_n)\) is not the
component of \(R-D\) whose root is at \(u_1\).  \end{lemma}

\begin{proof}  The proof is dual to that of the previous lemma.  If
\(V\) is the component of \(R-D\) whose root is at \(u_1\) and \(U\)
is the component of \(R-D\) that contains the leaf \(\sigma(u_n)\)
of \(R\), then an iterated augmentation by \(V\) along \(u_1,
\ldots, u_n\) would leave both the imbalance and the number of
components of \(D-R\) unchanged, and would reduce the number of
components of \(R-D\).  \end{proof}

\begin{lemma} For each non-trivial component \(U\) of \(D-R\) there
is a unique leaf \(\lambda(U)\) of \(U\) so that if \(r(U)\) is the
root of \(U\), then there is an iterated augmentation chain
\(\lambda(U) = u_1, \ldots, u_n\) with \(\sigma(u_n)=r(U)\).

For each non-trivial component \(V\) of \(R-D\) there is a unique
leaf \(\lambda(V)\) of \(V\) so that if \(r(V)\) is the root of
\(V\), then there is an iterated augmentation chain \(r(V) = u_1,
\ldots, u_n\) with \(\sigma(u_n)=\lambda(V)\).  \end{lemma}

\begin{proof} We consider the first claim.  The argument for the
second is obtained from the first by making a few mechanical
substitutions.

Consider \(u_1=\sigma^{-1}(r(U))\).  From previous lemmas, \(u_1\)
can only be a leaf of \(D\) that is also a leaf of \(R\), or it is a
leaf of \(U\).

Now, assume that we have an iterated augmentation chain \(u_1,
\ldots, u_n\) with \(r(U)=\sigma(u_n)\) and with \(u_1\) either a
leaf of \(U\) or a leaf of \(R\).  From the previous paragraph, we
know that there is at least one such chain.  In the case that
\(u_1\) is a leaf of \(U\), we are done.

Assume that \(u_1\) is a leaf of \(R\).  Increase the subscript of
every element of the chain by 1, so the chain now starts with
\(u_2\) and consider \(u_1=\sigma^{-1}(u_2)\).  From previous
lemmas, \(u_1\) can only be a leaf of \(D\) that is also a leaf of
\(R\), or it is a leaf of \(U\).  Thus we have created a chain that
is one longer than the original.

Since the element represented by \((D, \sigma, R)\) is a
homeomorphism, the chain cannot extend to an infinite loop.  Since
there are finitely many leaves of \(D\), this process must stop.
\end{proof}

We now make some definitions.  A leaf \(u\) of \(D\) is called
\begin{enumerate} \item {\itshape neutral} if it is a leaf of \(R\),
\item a {\itshape repeller} if it is \(\lambda(U)\) of some
component \(U\) of \(D-R\), and \item a {\itshape source} if it is a
leaf other than \(\lambda(U)\) of some component \(U\) of \(D-R\),
\item a {\itshape domain of attraction} if it is \(r(V)\) of some
component \(V\) of \(R-D\).  \end{enumerate} From what we have
shown, these cases are exhaustive and mutually exclusive.

A leaf \(v\) of \(R\) is called \begin{enumerate} \item {\itshape
neutral} if it is a leaf of \(D\), \item an {\itshape attractor} if
it is \(\lambda(V)\) of some component \(V\) of \(R-D\), and \item a
{\itshape sink} if it is a leaf other than \(\lambda(V)\) of some
component \(V\) of \(R-D\), \item a {\itshape range of repulsion} if
it is \(r(U)\) of some component \(U\) of \(D-R\).  \end{enumerate}
From what we have shown, these cases are exhaustive and mutually
exclusive.

\begin{lemma}\mylabel{LastTechLemma} The imbalance is the number of
sinks and is also the number of sources.  \end{lemma}

\begin{proof} The imbalance is the number of carets summed over all
components of \(D-R\).  If a component \(U\) of \(D-R\) has \(n\)
carets, then it has \(n+1\) leaves, of which exactly one is
\(\lambda(U)\) and \(n\) are sources.  Thus the imbalance is the
number of sources.  Similarly, it is the number of sinks.
\end{proof}

\begin{proof}[Proof of Proposition \ref{BddOrbits}] We are now ready
to describe the dynamics of an \(f\) whose representative \((D,
\sigma, R)\) has been chosen as described before Lemma
\ref{FirstTechLemma}.  From Lemmas \ref{FirstTechLemma} through
\ref{LastTechLemma}, we know the following since \(\sigma\) is a
bijection from the leaves of \(D\) to the leaves of \(R\).

The first non-neutral leaf in the forward orbit of a repeller
\(\lambda(U)\) of \(D\) is the range of repulstion \(r(U)\) of \(R\).  
The first non-neutral leaf in the forward orbit of a source of \(D\)
is a sink of \(R\).
The first non-neutral leaf in the forward orbit of a domain of
attraction \(r(V)\) of \(D\) is the attractor \(\lambda(V)\) of \(R\).
A neutral leaf of \(D\) that is not in the forward orbit of a
repeller, source or domain of attraction is in a finite cyclic
orbit of neutral leaves. 

The interval \(I\) in the Cantor set represented by a repeller
\(\lambda(U)\) is properly contained in the interval \(J\)
represented by \(r(U)\) and is taken affinely to it by an iterate of
\(f\).  Thus \(I\) contains a unique periodic point of \(f\) which
is represented by \(r(U)\) followed by infinite repetitions of the
path from \(r(U)\) to \(\lambda(U)\).  The period of the periodic
point is the number of iterations of \(f\) required to take
\(\lambda(U)\) to \(r(U)\).

The interval \(I\) in the Cantor set represented by a domain of
attraction \(r(V)\) properly contains the interval \(J\) represented
by \(\lambda(V)\) and is taken affinely to it by an iterate of
\(f\).  Thus \(I\) contains a unique periodic point of \(f\) which
is represented by \(r(V)\) followed by infinite repetitions of the
path from \(r(V)\) to \(\lambda(V)\).  The period of the periodic
point is the number of iterations of \(f\) required to take \(r(V)\)
to \(\lambda(V)\).

A sink of \(R\) is a leaf other than \(\lambda(V)\) of a component
\(V\) of \(R-D\), and the points in the interval \(I\) corresponding
to the sink are contained in the larger interval \(J\) corresponding
to \(r(V)\) of \(D\).  A finite iteration of \(f\) takes all of
\(J\) into the interval \(K\) corresponding to the attractor
\(\lambda(V)\) of \(R\).  Since \(I\) and \(K\) are disjoint, no
point in \(I\) is periodic.  Since every source has a sink in its
forward orbit, no point in an interval corresponding to a source of
\(D\) is periodic.

The neutral leaves not involved in any of the above will be
organized into a finite number of circuits, each of finite period.

The number of periods in the above discussion is finite which
completes the proof.  \end{proof}

\section{More on the baker's map}

In the next section, we give a simple proof of the simplicity of
\(V\).  It is much simpler than our proof the simplicity of \(2V\).
As an obstruction to an analogous proof for \(2V\) stands the
baker's map.

Since \(2V\) is simple, the baker's map must be a product of
commutators.  Before doing the calculations showing that \(2V\)
equals its commutator subgroup, the author did not believe that it
was.  A calculation that is an extension of the proof of Proposition
\ref{TwoVeeAbel} yields the following expression of \(C_0\) (the
baker's map) as a product of commutators:

\mymargin{BakerComm}\begin{equation}\label{BakerComm}
\begin{split}
C_0 =&\opi_1(K_6K_1K_2 K_3 K_4 K_5K_1K_2K_8
K_4^{-1})\opi_1^{\,-1}
K_2^{-1}K_1^{-1}K_6K_6K_1K_2 \\
&\quad K_3K_4 K_1K_2K_8K_7 
(K_1K_2 K_3 K_4 K_5K_1K_2)^{-1}.
\end{split}
\end{equation}

In that expression 
\[\begin{split}
K_1 &= [A_0^{-1}, A_1], \\
K_2 &= A_1[\pi_1^{-1}, A_0^{-1}]A_1^{-1}, \\
K_3 &= [\pi_1^{-1}, \opi_1^{\,-1}], \\
K_4 &= [\opi_1^{\,-1}, A_0^{-1}],
\end{split}
\qquad
\begin{split}
K_5 &= [A_1^{-1},\pi_0^{-1}], \\
K_6 &= [A_0^{-1},C_1^{-1}], \\
K_7 &= [B_1^{-1}, \pi_0^{-1}], \\
K_8 &= [\pi_1^{-1},A_0^{-1}].
\end{split}\]

The author has no idea if there is a shorter word than
\tref{BakerComm} expressing the baker's map as a product of
commutators.  The author found \tref{BakerComm} quite useless as an
aid to understanding why the baker's map was in the commutator
subgroup.

\section{The simplicity of \protect\(V\protect\) and the problem of
the baker's map}

We give a short proof of the simplicity of \(V\).  It is along the
lines of the proof in \cite{brown:finiteprop} but seems somewhat
shorter.  It is much shorter than the proof in \cite{CFP} where the
proof of simplicity is intimately tied to the calculation of a
finite presentation for \(V\).  The proof below is due to Mati Rubin
and is included here with his permission.  The basic idea is that
\(V\) is mostly a permutation group and all the permutations can be
made even by augmentation.

We say that an element of \(V\) is a \emph{permutation} if it is
represented by a triple of the form \((D, \sigma , D)\).  We say
that a caret \((u, u0, u1)\) in a finite tree \(D\) is
\emph{exposed} if both \(u0\) and \(u1\) are leaves of \(D\).  It is
clear that every finite tree has at least one exposed caret.

\begin{lemma}\mylabel{PermsGenVee}  The permutations generate
\(V\).  \end{lemma}

\begin{proof} Let \((D, \sigma, R)\) represent an arbitrary element
\(f\) of \(V\).  There is a permutation \((D, \tau, D)\) so that
\((D, \tau, D)(D, \sigma, R)= (D, \sigma', R)\) in which \(\sigma'\)
takes the leaves in an exposed caret of \(D\) in left-right order to
the leaves in an exposed caret of \(R\).  Now the reverse of an
augmentation can remove this exposed caret from both \(D\) and \(R\)
and we are done by induction since we will ultimately reduce \(D\)
and \(R\) to a tree with no carets.  \end{proof}

We say that a non-identity permutation \((D, \sigma, D)\) is a
\emph{transposition} if \(\sigma\) fixes all but two leaves of
\(D\).  The transposition is \emph{proper} if \(D\) has at least
three leaves.

\begin{lemma}\mylabel{PropTransGenVee} The proper transpositions
generate \(V\).  \end{lemma}

\begin{proof} By augmenting, every permutation can be represented as
\((D, \sigma, D)\) where \(D\) has at least three leaves.  The
result is clear.  \end{proof}

It is immediate that all of the proper transpositions are
conjugate.  Thus a normal subgroup of \(V\) is all of \(V\) if it
contains a proper transposition.  We now show that the normal
closure of any non-trivial element contains a proper transposition.

A non-trivial \(f\) with normal closure \(N\) moves intervals in
patterns, but we can think of this as moving nodes in the infinite
binary tree \(T\).  Since \(f\) must move some interval off itself,
it must move some caret \((u, u0, u1)\) to some caret \((v, v0,
v1)\) so that neither \(u\) nor \(v\) is in a subtree of \(T\) with
the other as root.  These can be chosen so that \(u\) and \(v\) are
at least distance three from the root of \(T\).  Let \(g\) be a
transposition interchanging \(u0\) with \(u1\).  Now \(h=[g,f]\)
interchanges \(u0\) with \(u1\) and \(v0\) with \(v1\).  Let \(j\)
be a transposition interchanging \(u0\) with \(v0\).  Now
\(k=[j,h]\) interchanges \(u0\) with \(v0\) and \(u1\) with \(v1\).
But this just interchanges \(u\) with \(v\) and is a proper
transposition since \(u\) and \(v\) are far enough from the root.
However, \(h\) and thus \(k\) are in \(N\).  We have finished the
proof of the following.

\begin{prop} The group \(V\) is simple.  \end{prop}

The outline above breaks down at the first step for \(2V\).  The
presence of the baker's map invalidates the proof of Lemma
\ref{PermsGenVee} for \(2V\).  It would be interesting to know if
there is a proof of the simplicity of \(2V\) that is shorter than
the one in this paper.

%\bibliography{thompson}

\providecommand{\bysame}{\leavevmode\hbox to3em{\hrulefill}\thinspace}

\noindent Department of Mathematical Sciences

\noindent State University of New York at Binghamton

\noindent Binghamton, NY 13902-6000

\noindent USA

\noindent email: matt@math.binghamton.edu

\end{document}